\documentclass{lms}
\usepackage{amssymb}
\usepackage{amsmath}
\usepackage[matrix,arrow,curve,cmtip]{xy}
\usepackage{mathrsfs}
\usepackage{graphicx}

\numberwithin{equation}{section}
\newtheorem{lettertheorem}{Theorem}

\newtheorem{theorem}[equation]{Theorem} 
\newtheorem{lemma}[equation]{Lemma}     
\newtheorem{cor}[equation]{Corollary}
\newtheorem{prop}[equation]{Proposition}

\newnumbered{definition}[equation]{Definition}
\newnumbered{remark}[equation]{Remark}
\newnumbered{example}[equation]{Example}

\title[$K$-theory of truncated polynomial algebras]{The tower of
$K$-theory of truncated polynomial algebras}

\author{Lars Hesselholt}

\classno{19D55 (primary), 19E15 (secondary)}

\extraline{The author was supported in part by the National Science
Foundation}

\begin{document}

\maketitle

\section*{Introduction}\label{introduction}

Suppose that $A$ is a regular noetherian ring that is also an
$\mathbb{F}_p$-algebra. Then it was proved by the author and
Madsen~\cite{hm1,hm2,h} that there is a long-exact sequence
$$\xymatrix{
{ \cdots } \ar[r] &
{ \bigoplus_{i \geqslant 0} \mathbb{W}_{i+1}\Omega_A^{q-2i} }
\ar[r]^(.48){V_m} &
{ \bigoplus_{i \geqslant 0} \mathbb{W}_{m(i+1)} \Omega_A^{q-2i} }
\ar[r]^(.48){\varepsilon} &
{ K_{q+1}(A[x]/(x^m),(x)) } \ar[r] &
{ \cdots } \cr
}$$
which expresses the $K$-groups of the truncated polynomial algebra
$A[x]/(x^m)$ relative to the ideal $(x)$ in terms of the groups
$\mathbb{W}_r\Omega_A^q$ of big de~Rham-Witt forms; see
also~\cite{h3}. We remark that in the papers~\cite{hm1,hm2,h} the
result was stated and proved under the more restrictive assumption
that the ring $A$ be a smooth $\mathbb{F}_p$-algebra. However, since
Popescu~\cite{popescu} has proved that every ring $A$ as above can be
written as a filtered colimit of smooth $\mathbb{F}_p$-algebras, the
general case immediately follows. In this paper, we consider the map
of relative $K$-groups
$$f_* \colon K_{q+1}(A[x]/(x^m),(x)) \to K_{q+1}(A[x]/(x^n),(x))$$
induced by the canonical projection
$$f \colon A[x]/(x^m) \to A[x]/(x^n).$$
To state the result, we recall from~\cite{hm2} that the big
de~Rham-Witt groups $\mathbb{W}_r\Omega_A^q$ are modules over the ring
$\mathbb{W}(A)$ of big Witt vectors in $A$. In particular, they are
modules over $\mathbb{W}(\mathbb{F}_p)$. The ring
$\mathbb{W}(\mathbb{F}_p)$ is canonically isomorphic to the product
indexed by the set of positive integers not divisible by $p$ of copies
of the ring $\mathbb{Z}_p$ of $p$-adic integers. A unit $\alpha$ of
the total quotient ring of $\mathbb{W}(\mathbb{F}_p)$ determines a
divisor $\operatorname{div}(\alpha)$ on $\mathbb{W}(\mathbb{F}_p)$
and, conversely, the unit $\alpha$ is determined, up to multiplication
by a unit of $\mathbb{W}(\mathbb{F}_p)$, by the divisor
$\operatorname{div}(\alpha)$. The ring $\mathbb{W}_r(\mathbb{F}_p)$ of
big Witt vectors of length $r$ in $\mathbb{F}_p$ determines a divisor
on $\mathbb{W}(\mathbb{F}_p)$ that we denote by
$\operatorname{div}(\mathbb{W}_r(\mathbb{F}_p))$.

\begin{lettertheorem}\label{main}Let $A$ be a regular noetherian ring
and an $\mathbb{F}_p$-algebra. Then the canonical projection $f \colon
A[x]/(x^m) \to A[x]/(x^n)$ induces a map of long-exact sequences
$$\xymatrix{
{ \cdots } \ar[r] &
{ \bigoplus_{i \geqslant 0} \mathbb{W}_{i+1}\Omega_A^{q-2i} } 
\ar[r]^(.48){V_m} \ar[d] &
{ \bigoplus_{i \geqslant 0} \mathbb{W}_{m(i+1)} \Omega_A^{q-2i} }
\ar[r]^(.49){\varepsilon} \ar[d] &
{ K_{q+1}(A[x]/(x^m),(x)) } \ar[r] \ar[d] &
{ \cdots } \cr
{ \cdots } \ar[r] &
{ \bigoplus_{i \geqslant 0} \mathbb{W}_{i+1}\Omega_A^{q-2i} } 
\ar[r]^(.49){V_n} &
{ \bigoplus_{i \geqslant 0} \mathbb{W}_{n(i+1)} \Omega_A^{q-2i} }
\ar[r]^(.49){\varepsilon} &
{ K_{q+1}(A[x]/(x^n),(x)) } \ar[r] &
{ \cdots } \cr
}$$
where the right-hand vertical map is the map of relative $K$-groups
induced by the canonical projection, where the middle vertical map
takes the $i$th summand of the domain to the $i$th summand of the
target by the composition
$$\xymatrix{
{ \mathbb{W}_{m(i+1)}\Omega_A^{q-2i} } \ar[r]^{\operatorname{res}} &
{ \mathbb{W}_{n(i+1)}\Omega_A^{q-2i} } \ar[r]^{m_{\alpha}} &
{ \mathbb{W}_{n(i+1)}\Omega_A^{q-2i} } \cr
}$$
of the restriction map and the multiplication by an element $\alpha =
\alpha_p(m,n,i)$ of $\mathbb{W}(\mathbb{F}_p)$ that is determined, up
to a unit, by the effective divisor
$$\operatorname{div}(\alpha) = \sum_{0 \leqslant h < i} \big(
\operatorname{div}(\mathbb{W}_{m(h+1)}(\mathbb{F}_p)) - 
\operatorname{div}(\mathbb{W}_{n(h+1)}(\mathbb{F}_p)) \big),$$
and where the left-hand vertical map is zero.
\end{lettertheorem}

The divisor $\operatorname{div}(\alpha_p(m,n,i))$ not only depends on
the integers $m$, $n$, and $i$, but also on the prime number
$p$. In more detail, recall that the length of a module $M$ over a
ring $R$ is the length $s$ of a longest ascending chain $M_0
\varsubsetneq M_1 \varsubsetneq \dots \varsubsetneq M_s$ of
$R$-submodules of $M$. If a longest chain does not exists, the module
$M$ is said to have infinite length. A $\mathbb{Z}_p$-module $M$ of
finite order $p^s$ has length $s$. Then the lengths as
$\mathbb{Z}_p$-module of the domain and target of the map
$$f_* \colon K_{2i+1}(\mathbb{F}_p[x]/(x^m),(x)) \to
K_{2i+1}(\mathbb{F}_p[x]/(x^n),(x))$$
are equal to $(m-1)(i+1)$ and $(n-1)(i+1)$, respectively, and hence,
do not depend on the prime $p$. By contrast, the lengths as a
$\mathbb{Z}_p$-module of the kernel
$K_{2i+1}(\mathbb{F}_p[x]/(x^m),(x^n))$ and cokernel
$K_{2i}(\mathbb{F}_p[x]/(x^m),(x^n))$ of this map are arithmetic
functions of the prime number $p$.

We note that, if $\operatorname{div}(\alpha_p(m,n,i)) \geqslant
\operatorname{div}(\mathbb{W}_{n(i+1)}(\mathbb{F}_p))$ then the
multiplication by $\alpha$ map
$$m_{\alpha} \colon \mathbb{W}_{n(i+1)}\Omega_A^{q - 2i} \to
\mathbb{W}_{n(i+1)}\Omega_A^{q-2i}$$
is zero. We examine the divisor $\alpha_p(m,n,i)$ and prove the
following result.

\begin{lettertheorem}\label{maink}Let $A$ be a regular noetherian ring
and an $\mathbb{F}_p$-algebra and assume that $A$ is finitely
generated as an algebra over the subring $A^p \subset A$ of $p$th
powers. Let $m$ and $n$ be positive integers with $m > n+1$. Then
there exists an integer $q_0 \geqslant 1$ such that the map
$$f_* \colon K_{q+1}(A[x]/(x^m),(x)) \to K_{q+1}(A[x]/(x^n),(x))$$
induced by the canonical projection is zero, for all
$q \geqslant q_0$.
\end{lettertheorem}

The main purpose of the paper~\cite{hm2} was to evaluate the
$\operatorname{Nil}$-groups of the truncated polynomial rings
$A[x]/(x^m)$. We recall that, for every ring $R$, the ring
homomorphisms $\eta \colon R \to R[t]$ and $\epsilon \colon R[t] \to
R$ defined by $\eta(a) = a$ and $\epsilon(f) = f(0)$ give rise to a
direct sum decomposition
$$K_{q+1}(R[t]) = K_{q+1}(R) \oplus K_{q+1}(R[t],(t)).$$
The second summand on the right-hand side is the
$\operatorname{Nil}$-group
$$\operatorname{Nil}_q(R) = K_{q+1}(R[t],(t))$$
which measures the extent to which the functor $K_{q+1}$ fails to be
homotopy invariant. The group $\operatorname{Nil}_q(R)$ is zero, if
the ring $R$ is regular. This implies that, for $R = A[x]/(x^m)$ with
$A$ regular, the canonical map
$$\operatorname{Nil}_q(A[x]/(x^m),(x)) \to
\operatorname{Nil}_q(A[x]/(x^m))$$
is an isomorphism. The ring homomorphisms $\eta \colon R \to R[t]$ and
$\epsilon \colon R[t] \to R$ also give rise to a direct sum
decomposition of big de~Rham-Witt groups 
$$\mathbb{W}_r\Omega_{R[t]}^q = \mathbb{W}_r\Omega_R^q \oplus
\mathbb{W}_r\Omega_{(R[t],(t))}^q.$$
Hence, if $A$ is a regular noetherian ring and an
$\mathbb{F}_p$-algebra, the long-exact sequence at the beginning of
the paper implies the long-exact sequence
$$\xymatrix{
{ \cdots } \ar[r] &
{ \bigoplus_{i \geqslant 0} \mathbb{W}_{i+1}\Omega_{(A[t],(t))}^{q-2i} }
\ar[r]^(.48){V_m} &
{ \bigoplus_{i \geqslant 0} \mathbb{W}_{m(i+1)} \Omega_{(A[t],(t))}^{q-2i} }
\ar[r]^(.57){\varepsilon} &
{ \operatorname{Nil}_q(A[x]/(x^m)) } \ar[r] &
{ \cdots } \cr
}$$
which expresses the $\operatorname{Nil}$-groups of $A[x]/(x^m)$ in
terms of the relative groups of big de~Rham-Witt forms. The map $V_m$
is injective, if $p$ does not divide $m$, but is generally not
injective, if $p$ divides $m$. We remark that~\cite[Thm.~B]{hm3} gives
an explicit description of the middle and left-hand terms of this
sequence as functors of the big de~Rham-Witt groups of the ring
$A$. In particular, the group $\operatorname{Nil}_q(A[x]/(x^m))$ is
non-zero, for all integers $q \geqslant 0$ and $m > 1$. We note that
$\operatorname{Nil}_q(A[x]/(x^m))$ is zero, for $q < 
0$~\cite[Chap.~XII, Prop.~10.1]{bass}. Thm.~\ref{main} above
immediately implies the following result.

\begin{lettertheorem}\label{mainnil}Let $A$ be a regular noetherian
ring and an $\mathbb{F}_p$-algebra. Then the canonical projection
$f \colon A[x]/(x^m) \to A[x]/(x^n)$ induces a map of long-exact
sequences
$$\xymatrix{
{ \cdots } \ar[r] &
{ \bigoplus_{i \geqslant 0} \mathbb{W}_{i+1}\Omega_{(A[t],(t))}^{q-2i} }  
\ar[r]^(.48){V_m} \ar[d] &
{ \bigoplus_{i \geqslant 0} \mathbb{W}_{m(i+1)} \Omega_{(A[t],(t))}^{q-2i} }
\ar[r]^(.57){\varepsilon} \ar[d] &
{ \operatorname{Nil}_q(A[x]/(x^m)) } \ar[r] \ar[d]^{f_*} &
{ \cdots } \cr
{ \cdots } \ar[r] &
{ \bigoplus_{i \geqslant 0} \mathbb{W}_{i+1}\Omega_{(A[t],(t))}^{q-2i} } 
\ar[r]^(.49){V_n} &
{ \bigoplus_{i \geqslant 0} \mathbb{W}_{n(i+1)} \Omega_{(A[t],(t))}^{q-2i} }
\ar[r]^(.57){\varepsilon} &
{ \operatorname{Nil}_q(A[x]/(x^n)) } \ar[r] &
{ \cdots } \cr
}$$
where the middle vertical map takes the $i$th summand of the domain to the
$i$th summand of the target by the composition
$$\xymatrix{
{ \mathbb{W}_{m(i+1)}\Omega_{(A[t],(t))}^{q-2i} } \ar[r]^{\operatorname{res}} &
{ \mathbb{W}_{n(i+1)}\Omega_{(A[t],(t))}^{q-2i} } \ar[r]^{m_{\alpha}} &
{ \mathbb{W}_{n(i+1)}\Omega_{(A[t],(t))}^{q-2i} } \cr
}$$
of the restriction map and the multiplication by the element $\alpha =
\alpha_p(m,n,i)$ of Thm.~\ref{main}, and where the left-hand vertical
map is zero.
\end{lettertheorem}

Similarly, Thm.~\ref{maink} above immediately implies the following
result.

\begin{lettertheorem}\label{mainknil}Let $A$ be a regular noetherian
ring and an $\mathbb{F}_p$-algebra and assume that $A$ is finitely
generated as an algebra over the subring $A^p \subset A$ of $p$th
powers. Let $m$ and $n$ be positive integers with $m > n+1$. Then
there exists an integer $q_0 \geqslant 1$ such that the map
$$f_* \colon \operatorname{Nil}_q(A[x]/(x^m)) \to
\operatorname{Nil}_q(A[x]/(x^n))$$
induced by the canonical projection is zero, for all
$q \geqslant q_0$.
\end{lettertheorem}

We remark that, if $A$ is an $\mathbb{F}_p$-algebra and $m$ is a power
of $p$, then $A[x]/(x^m)$ is equal to the group algebra $A[C_m]$ of
the cyclic group of order $m$ with generator $1+x$.

All rings considered in this paper are assumed to be commutative and
unital. We write $\mathbb{T}$ for the multiplicative group of complex
numbers of modulus $1$. 

Finally, the author would like to express his gratitude to an
anonymous referee for a number of helpful suggestions on improving the
exposition of the paper.

\section{$p$-typical decompositions}\label{ptypical}

The long-exact sequences that appear in Thm.~\ref{main} of the
introduction have canonical $p$-typical product decompositions.
We here recall the $p$-typical decomposition and state the equivalent
Thm.~\ref{ptypicalmain}.

We recall from~\cite[Cor.~1.2.6]{hm2} that the big de~Rham-Witt groups
$\mathbb{W}_r\Omega_A^q$ decompose as a product of the $p$-typical
de~Rham-Witt groups $W_s\Omega_A^q$ of
Bloch-Deligne-Illusie~\cite{illusie}. More generally, there is a big
de~Rham-Witt group $\mathbb{W}_S\Omega_A^q$ associated to every subset
$S \subset \mathbb{N}$ of the set of positive integers stable under
division. The group $\mathbb{W}_{\emptyset}\Omega_A^q$ is zero, the
group $\mathbb{W}_r\Omega_A^q$ is the big de~Rham-Witt group
associated to the set of positive integers less than or equal to
$r$, and the $p$-typical de~Rham-Witt group $W_s\Omega_A^q$ is the 
big de~Rham-Witt group
$$W_s\Omega_A^q = \mathbb{W}_{\{1,p,\dots,p^{s-1}\}}\Omega_A^q.$$
For every pair of subset $T \subset S \subset \mathbb{N}$ stable under
division, there is a map
$$\operatorname{res} \colon \mathbb{W}_S\Omega_A^q \to
\mathbb{W}_T\Omega_A^q \hskip7mm \text{(restriction)}.$$
We mention that, if $S \subset \mathbb{N}$ is the union of a family of
subsets $S_{\alpha} \subset \mathbb{N}$ each of which is stable under
division, then the restriction maps induce an isomorphism
$$\mathbb{W}_S\Omega_A^q \xrightarrow{\sim} \lim_{\operatorname{res}}
\mathbb{W}_{S_{\alpha}}\Omega_A^q.$$
We also mention that every subset $S \subset \mathbb{N}$ stable under
division is equal to the union of the subsets $\langle s \rangle
\subset \mathbb{N}$, $s \in S$, where $\langle s \rangle$ is the set
of divisors of $s$. Let $S \subset \mathbb{N}$ be a subset stable
under division, let $s$ be any positive integer, and let $S/s$ be the
set of positive integers $t$ such that $st \in S$. Then there are maps
$$\begin{aligned}
F_s \colon & \mathbb{W}_S\Omega_A^q \to \mathbb{W}_{S/s}\Omega_A^q
\hskip7mm \text{(Frobenius)} \cr
V_s \colon & \mathbb{W}_{S/s}\Omega_A^q \to \mathbb{W}_S\Omega_A^q
\hskip7mm \text{(Verschiebung).} \cr
\end{aligned}$$
Now, let $S \subset \mathbb{N}$ be a finite subset stable under
division, and let $j$ be a positive integer that is not divisible by
$p$. We define the non-negative integer $u = u_p(S,j)$ by
$$u_p(S,j) = 
\operatorname{card}( S/j \cap \{1, p, p^2, \dots \} ) =
\operatorname{card}( S \cap \{j, pj, p^2j, \dots \} )$$
and note that $u$ is non-zero if and only if $j \in S$. We then let
$$\eta_j \colon \mathbb{W}_S\Omega_A^q \to W_u\Omega_A^q$$
be the composite map
$$\xymatrix{
{ \mathbb{W}_S\Omega_A^q } \ar[r]^(.47){F_j} &
{ \mathbb{W}_{S/j}\Omega_A^q } \ar[r]^(.39){\operatorname{res}} &
{ \mathbb{W}_{\{1,p,\dots,p^{u-1}\}}\Omega_A^q } \ar@{=}[r] &
{ W_u\Omega_A^q } \cr
}$$
and define
$$\eta \colon \mathbb{W}_S\Omega_A^q \to \prod_j W_u\Omega_A^q$$
to be the map given by the maps $\eta_j$ as $j$ ranges over the set
$I_p$ of all positive integers not divisible by $p$. The map $\eta$ is
an isomorphism, for every $\mathbb{Z}_{(p)}$-algebra,
by~\cite[Cor.~1.2.6]{hm2}. We remark that, on the right-hand side,
the factors indexed by $j \notin S$ are zero. The maps
$\operatorname{res}$, $F_s$, and $V_s$ are also expressed as products
of their $p$-typical analogs
$$\begin{aligned}
R = \operatorname{res} \colon & W_s\Omega_A^q \to W_{s-1}\Omega_A^q
\hskip6.7mm \text{(restriction)} \hfill\space \cr
F = F_p \colon & W_s\Omega_A^q \to W_{s-1}\Omega_A^q
\hskip7mm \text{(Frobenius)} \hfill\space \cr
V = V_p \colon & W_{s-1}\Omega_A^q \to W_s\Omega_A^q
\hskip7mm \text{(Verschiebung).} \hfill\space \cr
\end{aligned}$$
Let $T \subset S \subset \mathbb{N}$ be a pair of subsets stable under
division. Then there is a commutative diagram
$$\xymatrix{
{ \mathbb{W}_S\Omega_A^q } \ar[r]^(.42){\eta}
\ar[d]^{\operatorname{res}} &
{ \prod_{j} W_u\Omega_A^q } \ar[d]^{\operatorname{res}^{\eta}} \cr
{ \mathbb{W}_T\Omega_A^q } \ar[r]^(.42){\eta} &
{ \prod_{j} W_{u'}\Omega_A^q, } \cr
}$$
where $u = u_p(S,j)$ and $u' = u_p(T,j)$, and where the map
$\operatorname{res}^{\eta}$ takes the factor indexed by an integer
$j \in T$ that is not divisible by $p$ to the factor indexed by the
same integer $j$ by the map $R^{u-u'}$ and annihilates the remaining
factors. Similarly, let $S \subset \mathbb{N}$ be a subset stable
under division, and let $s$ be a positive integer. We write
$s = p^vs'$ with $s'$ not divisible by $p$. Then there are commutative
diagrams
$$\xymatrix{
{ \mathbb{W}_S\Omega_A^q } \ar[r]^(.42){\eta} \ar[d]<1ex>^{F_s} &
{ \prod_j W_u\Omega_A^q } \ar[d]<1ex>^{F_s^{\eta}} \cr
{ \mathbb{W}_{S/s}\Omega_A^q } \ar[r]^(.42){\eta} \ar[u]<1ex>^{V_s} &
{ \prod_j W_{u-v}\Omega_A^q, } \ar[u]<1ex>^{V_s^{\eta}} \cr
}$$
where $u = u_p(S,j)$, and where the maps $F_s^{\eta}$ and $V_s^{\eta}$
are defined as follows: The map $F_s^{\eta}$ takes the factor indexed
by an integer $j \in S$ that is not divisible by $p$, but is divisible
by $s'$, and satisfies $p^vj \in S$ to the factor indexed by the
integer $j/s' \in S/s$ by the map $F^v$ and annihilates the remaining
factors. The map $V_s^{\eta}$ takes the factor indexed by $j \in S/s$
not divisible by $p$ to the factor indexed by $s'j \in S$ by the map
$s'V^v$.

Specializing to the case $S = \{1,2,\dots,m(i+1)\}$, we define
\begin{equation}\label{s_p(m,i,j)}
s_p(m,i,j) = u_p(S,j) = \operatorname{card}( \{1,2,\dots,m(i+1)\}
\cap \{j, pj, p^2j, \dots \}).
\end{equation}
The integer $s_p(m,i,j)$ is non-zero if and only if $1 \leqslant j
\leqslant m(i+1)$, and in this case, is the unique integer $s$
such that $p^{s-1}j \leqslant m(i+1) < p^sj$.
Suppose now that $A$ is a regular noetherian ring and an
$\mathbb{F}_p$-algebra. It is then straightforward to verify that the
$p$-typical decomposition of the top long-exact sequence in the
diagram in Thm.~\ref{main} of the introduction takes the form
\begin{equation}\label{mdrw}
\begin{aligned}
\cdots & \to
\bigoplus_{i \geqslant 0} \bigoplus_{j \in m'I_p}
W_{s-v}\Omega_A^{q-2i} \xrightarrow{m'V^v} \bigoplus_{i \geqslant 0}
\bigoplus_{j \in I_p} W_s\Omega_A^{q-2i} \cr 
{} & \xrightarrow{\varepsilon} \;\,
{ K_{q+1}(A[x]/(x^m),(x)) } \;\, \xrightarrow{\partial} \;\;\,
\bigoplus_{i \geqslant 0} \bigoplus_{j \in m'I_p}
W_{s-v}\Omega_A^{q-1-2i} \to \cdots \cr
\end{aligned}
\end{equation}
where $s = s_p(m,i,j)$, where $m = p^vm'$ with $m'$ not divisible
by $p$, and where $I_p$ is the set of positive integers not divisible
by $p$. Similarly, the bottom long-exact sequence in the diagram in
the statement of Thm.~\ref{main} takes the form
\begin{equation}\label{ndrw}
\begin{aligned}
\cdots & \to
\bigoplus_{i \geqslant 0} \bigoplus_{j \in n'I_p}
W_{t-w}\Omega_A^{q-2i} \xrightarrow{n'V^w} \bigoplus_{i \geqslant 0}
\bigoplus_{j \in I_p} W_t\Omega_A^{q-2i} \cr 
{} & \xrightarrow{\varepsilon} \;\,
{ K_{q+1}(A[x]/(x^n),(x)) } \;\, \xrightarrow{\partial} \;\;\,
\bigoplus_{i \geqslant 0} \bigoplus_{j \in n'I_p}
W_{t-w}\Omega_A^{q-1-2i} \to \cdots \cr
\end{aligned}
\end{equation}
where $t = s_p(n,i,j)$, and where $n = p^wn'$ with $n'$ not divisible
by $p$. 

Finally, we explain the $p$-typical decomposition of the divisor
$\operatorname{div}(\alpha_p(m,n,i))$ that appears in
Thm.~\ref{main}. Let $\mathscr{S} \subset \mathbb{W}(\mathbb{F}_p)$ be
the subset of non-zero-divisors. Then the map
$$\operatorname{div} \colon
(\mathscr{S}^{-1}\mathbb{W}(\mathbb{F}_p))^* /
\mathbb{W}(\mathbb{F}_p)^* \to
\operatorname{Div}(\mathbb{W}(\mathbb{F}_p))$$
is injective. We refer to~\cite[\S 21]{ega4} for the general theory of
divisors. Let also $S \subset W(\mathbb{F}_p)$ be the subset of
non-zero-divisors, and let
$$\operatorname{ord} \colon
(\mathscr{S}^{-1}\mathbb{W}(\mathbb{F}_p))^* /
\mathbb{W}(\mathbb{F}_p)^* \xrightarrow{\sim}
\prod_{j \in I_p} (S^{-1}W(\mathbb{F}_p))^* / W(\mathbb{F}_p)^*
\xrightarrow{\sim}
\prod_{j \in I_p} \mathbb{Z}$$
be the composition of the isomorphism induced by the ring isomorphism
$\eta$ and the isomorphism given by the $p$-adic valuation. We recall
the function $s_p(m,i,j)$ from~(\ref{s_p(m,i,j)}) above. The
$p$-typical decomposition of the ring
$\mathbb{W}_{m(i+1)}(\mathbb{F}_p)$ immediately implies the following
result.

\begin{lemma}\label{ptypicaldivisor}For every $a \in 
(\mathscr{S}^{-1}\mathbb{W}(\mathbb{F}_p))^*$, the following are
equivalent:
\begin{enumerate}
\item[(i)] $\operatorname{div}(a) =
\operatorname{div}(\mathbb{W}_{m(i+1)}(\mathbb{F}_p))$.
\item[(ii)] $\operatorname{ord}(a) = (s_p(m,i,j) \mid j \in I_p)$.
\end{enumerate}
\end{lemma}

The following statement is equivalent to Thm.~\ref{main} of the
introduction. The proof is given in Sect.~\ref{proofoftheorem} below.

\begin{theorem}\label{ptypicalmain}Let $A$ be a regular noetherian
ring and an $\mathbb{F}_p$-algebra. Then the canonical projection
$f \colon A[x]/(x^m) \to A[x]/(x^n)$ induces a map of long-exact
sequences from the sequence~(\ref{mdrw}) to the sequence~(\ref{ndrw})
that is given, on the lower left-hand terms, by the map
$$f_* \colon K_*(A[x]/(x^m),(x)) \to K_*(A[x]/(x^n),(x))$$
induced by the canonical projection, on the upper right-hand terms, by
the map that takes the $(i,j)$th summand of the domain to the
$(i,j)$th summand of the target by the composite map
$$\xymatrix{
{ W_s\Omega_A^{q - 2i} } \ar[r]^{R^{s-t}} &
{ W_t\Omega_A^{q - 2i} } \ar[r]^{m_{\alpha}} &
{ W_t\Omega_A^{q - 2i}, } \cr
}$$
where the right-hand map is multiplication by a $p$-adic integer
$\alpha = \alpha_p(m,n,i,j)$ whose $p$-adic valuation is given by the
sum
$$v_p(\alpha_p(m,n,i,j)) = \sum_{0 \leqslant h < i} \big( s_p(m,h,j) -
s_p(n,h,j) \big),$$
and, on the upper left-hand terms, by the zero map. Here $s_p(m,i,j)$
is the integer~(\ref{s_p(m,i,j)}).
\end{theorem}

\begin{remark}(i) The integer $s_p(m,h,j) - s_p(n,h,j)$ is equal to
the number of positive integers $r$ such that $n(h+1) < p^{r-1}j
\leqslant m(h+1)$. 

(ii) We recall from~\cite[Prop.~I.3.4]{illusie} that the map
$m_{\alpha} \colon W_t\Omega_A^{q-2i} \to W_t\Omega_A^{q-2i}$ given by
the multiplication by a $p$-adic integer $\alpha$ of $p$-adic
valuation $v$ factors canonically
$$\xymatrix{
{ W_t\Omega_A^{q - 2i} } \ar[r]^(.48){R^v} &
{ W_{t-v}\Omega_A^{q - 2i} } \ar[r]^(.55){\underline{m}_{\alpha}} &
{ W_t\Omega_A^{q-2i} } \cr
}$$
as the composite of the surjective restriction map $R^v$ and an
injective map $\underline{m}_{\alpha}$. There is no analog of this
factorization for the big de~Rham-Witt groups. Indeed, the quotient
of $\mathbb{W}(\mathbb{F}_p)$ by the ideal generated by the element
$\alpha_p(m,n,i) \in \mathbb{W}(\mathbb{F}_p)$ that appears in
Thm.~\ref{main} is generally not of the form
$\mathbb{W}_S(\mathbb{F}_p)$ for some subset $S\subset \mathbb{N}$
that is stable under division. 
\end{remark}

\section{Topological Hochschild homology}\label{trsection}

The proof of Thm.~\ref{main} of the introduction is based on the
following result which we proved in~\cite[Prop.~4.2.3]{hm1}: For every
$\mathbb{F}_p$-algebra $A$, there is a long-exact sequence
$$\xymatrix{
{ \cdots } \ar[r] &
{ \lim_R\operatorname{TR}_{q-\lambda_d}^{r/m}(A) } \ar[r]^{V_m} &
{ \lim_R\operatorname{TR}_{q-\lambda_d}^r(A) }
\ar[r]^(.45){\varepsilon} &
{ K_{q+1}(A[x]/(x^m),(x)) } \ar[r] &
{ \cdots, } \cr
}$$
which expresses the groups $K_{q+1}(A[x]/(x^m),(x))$ in terms of the
$RO(\mathbb{T})$-graded equivariant homotopy groups of the topological
Hochschild $\mathbb{T}$-spectrum $T(A)$. (The corresponding long-exact
sequence of homotopy groups with finite coefficients is valid for
every ring $A$.) We briefly recall the terms in this sequence.

The topological Hochschild $\mathbb{T}$-spectrum $T(A)$ associated
with the ring $A$ is a cyclotomic spectrum in the sense
of~\cite[Def.~2.2]{hm}. In particular, it is an object of the
$\mathbb{T}$-stable homotopy category. Let $\lambda$ be a finite
dimensional orthogonal $\mathbb{T}$-representation, let $S^{\lambda}$
be the one-point compactification, and let $C_r \subset \mathbb{T}$ be
the subgroup of order $r$. Then one defines
$$\operatorname{TR}_{q-\lambda}^r(A) = [ S^q \wedge (\mathbb{T}/C_r)_+, 
S^{\lambda} \wedge T(A) ]_{\mathbb{T}}$$
to be the abelian group of maps in the $\mathbb{T}$-stable homotopy
category between the indicated $\mathbb{T}$-spectra. Suppose that
$r = st$. Then there are maps
$$\begin{aligned}
F_s \colon & \operatorname{TR}_{q-\lambda}^r(A) \to
\operatorname{TR}_{q-\lambda}^t(A) \hskip10.4mm \text{(Frobenius)}
\hfill\space \cr
V_s \colon & \operatorname{TR}_{q-\lambda}^t(A) \to
\operatorname{TR}_{q-\lambda}^r(A) \hskip10.4mm \text{(Verschiebung)}
\hfill\space \cr
\end{aligned}$$
induced by maps
$f_s \colon (\mathbb{T}/C_t)_+ \to (\mathbb{T}/C_r)_+$ and 
$v_s \colon (\mathbb{T}/C_r)_+ \to (\mathbb{T}/C_t)_+$ in the
$\mathbb{T}$-stable homotopy category that are defined as follows. The
map $f_s$ is the map of suspension $\mathbb{T}$-spectra induced by the
canonical projection  $\operatorname{pr} \colon \mathbb{T}/C_t \to
\mathbb{T}/C_r$. The map $v_s$ is the corresponding transfer map. To
define it, we choose an embedding $\iota \colon \mathbb{T}/C_t
\hookrightarrow \lambda$ into a finite dimensional orthogonal
$\mathbb{T}$-presentation. The product embedding
$(\iota,\operatorname{pr}) \colon \mathbb{T}/C_t \to \lambda \times
\mathbb{T}/C_r$ has trivial normal bundle, and the linear structure on
$\lambda$ determines a preferred trivialization. Hence, the
Pontryagin-Thom construction gives a map of pointed
$\mathbb{T}$-spaces
$$S^{\lambda} \wedge (\mathbb{T}/C_r)_+ \to S^{\lambda} \wedge
(\mathbb{T}/C_t)_+$$
and $v_s$ is the induced map of suspension $\mathbb{T}$-spectra.

The isomorphism $\rho_s \colon \mathbb{T} \to \mathbb{T}/C_s$ given
by the $s$th root induces an equivalence of categories $\rho_s^*$ from
the $\mathbb{T}/C_s$-stable homotopy category to the
$\mathbb{T}$-stable homotopy category. The additional cyclotomic
structure of the $\mathbb{T}$-spectrum $T(A)$ gives rise to a map of
$\mathbb{T}$-spectra
$$r_s \colon \rho_s^*(T(A)^{C_s}) \to T(A).$$
Suppose again that $r = st$. Then we get a map of
$RO(\mathbb{T})$-graded equivariant homotopy groups
$$R_s \colon \operatorname{TR}_{q-\lambda}^r(A) \to
\operatorname{TR}_{q-\lambda'}^t(A) \hskip9.6mm \text{(restriction)},$$
where $\lambda' = \rho_s^*(\lambda^{C_s})$ is the
$\mathbb{T}/C_s$-representation $\lambda^{C_s}$ considered as a
$\mathbb{T}$-representation via the isomorphism $\rho_s$. The map
$R_s$ is defined to be the composition
$$\begin{aligned}
{} & \operatorname{TR}_{q-\lambda}^r(A) =
[ S^q \wedge (\mathbb{T}/C_r)_+, 
S^{\lambda} \wedge T(A) ]_{\mathbb{T}} \xleftarrow{\sim}
[ S^q \wedge (\mathbb{T}/C_r)_+, 
S^{\lambda^{C_s}} \wedge T(A)^{C_s} ]_{\mathbb{T}/C_s} \cr
& {} \xrightarrow{\sim} [ S^q \wedge (\mathbb{T}/C_t)_+,
S^{\lambda'} \wedge \rho_s^*(T(A)^{C_s}) ]_{\mathbb{T}} \to
[ S^q \wedge (\mathbb{T}/C_t)_+,
S^{\lambda'} \wedge T(A) ]_{\mathbb{T}} =
\operatorname{TR}_{q-\lambda'}^t(A) \cr
\end{aligned}$$
of the canonical isomorphism, the isomorphism $\rho_s^*$, and the map
induced by $r_s$. It is proved in~\cite[Addendum~3.3]{hm} that there
is a canonical  isomorphism of rings
$$\xi \colon \mathbb{W}_{\langle r\rangle}(A) \xrightarrow{\sim}
\operatorname{TR}_0^r(A)$$
from the ring of big Witt vectors in $A$ corresponding to the subset
$\langle r \rangle \subset \mathbb{N}$ of divisors of $r$. The
isomorphism $\xi$ is compatible with the restriction, Frobenius, and
Verschiebung operators.

We return to the long-exact sequence of~\cite[Prop.~4.2.3]{hm1}
displayed at the beginning of the section. In the middle and left-hand
terms, the representation $\lambda_d$ is defined to be the sum
$$\lambda_d = \mathbb{C}(d) \oplus \mathbb{C}(d-1) \oplus \dots \oplus
\mathbb{C}(1),$$
where $\mathbb{C}(i) = \mathbb{C}$ with $\mathbb{T}$ acting from the
left by $z \cdot w = z^iw$, and where
$$d = d(m,r) = \left[ \frac{r-1}{m} \right]$$
is the largest integer less than or equal to $(r-1)/m$. The limits
range over the set of all positive integers and the set of positive
integers divisible by $m$, respectively, and the structure map in both
limit systems is the restriction map. We note that if $r = st$ then
$$\rho_s^*(\lambda_{d(m,r)}^{C_s}) = \lambda_{d(m,t)}$$
as required. To prove Thm.~\ref{main} of the introduction, we first
prove the following general result which does not require the ring $A$
to be regular or noetherian. The proof of this result is given in
Sect.~\ref{barconstruction} below.

\begin{theorem}\label{maintr}Let $A$ be an\, $\mathbb{F}_p$-algebra,
and let $f \colon A[x]/(x^m) \to A[x]/(x^n)$ be the canonical
projection. Then there is a map of long-exact sequences
$$\xymatrix{
{ \cdots } \ar[r] &
{ \lim_R\operatorname{TR}_{q-\lambda_d}^{r/m}(A) } \ar[r]^{V_m} \ar[d] &
{ \lim_R\operatorname{TR}_{q-\lambda_d}^r(A) }
\ar[r]^(.45){\varepsilon} \ar[d] &
{ K_{q+1}(A[x]/(x^m),(x)) } \ar[r] \ar[d] &
{ \cdots } \cr
{ \cdots } \ar[r] &
{ \lim_R\operatorname{TR}_{q-\lambda_e}^{r/n}(A) } \ar[r]^{V_n} &
{ \lim_R\operatorname{TR}_{q-\lambda_e}^r(A) }
\ar[r]^(.45){\varepsilon} & 
{ K_{q+1}(A[x]/(x^n),(x)) } \ar[r] &
{ \cdots, } \cr
}$$
where $d = d(m,r) = [(r-1)/m]$ and $e = d(n,r) = [(r-1)/n]$, where the
right-hand vertical map is the map of relative $K$-groups induced by
the canonical projection, where the middle vertical map is the map of
limits defined by the maps
$$\iota(m,n,r)_q \colon \operatorname{TR}_{q-\lambda_d}^r(A) \to
\operatorname{TR}_{q-\lambda_e}^r(A)$$
induced by the canonical inclusions $\iota(m,n,r) \colon S^{\lambda_d}
\to S^{\lambda_e}$, and where the left-hand vertical map is zero. The
corresponding statement for the homotopy groups with finite
coefficients is valid for every ring $A$.
\end{theorem}

We show in Lemma~\ref{limitreached} below that the limits that appear
in the statement of Thm.~\ref{maintr} stabilize. In particular, the
corresponding derived limits vanish.

The long-exact sequences of Thm.~\ref{maintr} admit a $p$-typical
decomposition analogous to the $p$-typical decomposition of the big
de~Rham-Witt groups which we discussed in Sect.~\ref{ptypical}
above. In the remainder of this section, we recall this decomposition
and state the equivalent version Thm.~\ref{ptypicalmaintr} of
Thm.~\ref{maintr}. 

It is proved in~\cite[Prop.~4.2.5]{hm1} that the groups
$\operatorname{TR}_{q-\lambda}^r(A)$, which appear in the statement of
Thm.~\ref{maintr}, decompose as a product of the $p$-typical groups
$$\operatorname{TR}_{q-\lambda}^u(A;p) =
\operatorname{TR}_{q-\lambda}^{p^{u-1}}(A) = 
[ S^q \wedge (\mathbb{T}/C_{p^{u-1}})_+, S^{\lambda} \wedge T(A)
]_{\mathbb{T}}.$$
To state the result, we first write $r = p^{u-1}r'$, where $r'$ is not
divisible by $p$. Let $j$ be a divisor of $r'$, and let $\lambda' =
\lambda'(j) = \rho_{r'/j}^*(\lambda^{C_{r'/j}})$. We let 
$$\gamma_j \colon \operatorname{TR}_{q-\lambda}^r(A) \to
\operatorname{TR}_{q-\lambda'}^u(A;p)$$ 
be the composite map
$$\xymatrix{
{ \operatorname{TR}_{q-\lambda}^r(A) } \ar[r]^{F_j} &
{ \operatorname{TR}_{q-\lambda}^{r/j}(A) } \ar[r]^(.48){R_{r'/j}} &
{ \operatorname{TR}_{q-\lambda'}^{p^{u-1}}(A) } \ar@{=}[r] &
{ \operatorname{TR}_{q-\lambda'}^u(A;p) } \cr
}$$
and define
$$ \gamma \colon \operatorname{TR}_{q-\lambda}^r(A) \to
\prod_{j} \operatorname{TR}_{q-\lambda'}^u(A;p)$$
be the product of the maps $\gamma_j$ as $j$ ranges over the divisors
of $r'$. The map $\gamma$ is an isomorphism, for every
$\mathbb{Z}_{(p)}$-algebra $A$, by~\cite[Prop.~4.2.5]{hm1}. We remark
that this decomposition is analogous to the $p$-typical decomposition
of the big de~Rham-Witt groups $\mathbb{W}_{\langle r\rangle}\Omega_A^q$
that we discussed in Sect.~\ref{ptypical} above. Indeed, let $j$ be a
positive integer not divisible by $p$. Then the integer $u_p(\langle r
\rangle,j)$ is equal to $u = v_p(r) + 1$, if $j$ divides $r'$, and is
zero, otherwise. The maps $R_s$, $F_s$, and $V_s$ are similarly
expressed as products of their $p$-typical analogs
$$\begin{aligned}
R = R_p \colon & \operatorname{TR}_{q-\lambda}^u(A;p) \to
\operatorname{TR}_{q-\lambda'}^{u-1}(A;p) \hskip6.3mm
\text{(restriction)} \hfill\space \cr
F = F_p \colon & \operatorname{TR}_{q-\lambda}^u(A;p) \to
\operatorname{TR}_{q-\lambda}^{u-1}(A;p) 
\hskip7mm \text{(Frobenius)} \hfill\space \cr
V = V_p \colon & \operatorname{TR}_{q-\lambda}^{u-1}(A;p) \to
\operatorname{TR}_{q-\lambda}^u(A;p) 
\hskip7mm \text{(Verschiebung).} \hfill\space \cr
\end{aligned}$$
Suppose that $r = st$ and write $s = p^vs'$ and $t = p^{u-v-1}t'$
with $s'$ and $t'$ not divisible by $p$. Then there are commutative
diagrams
$$\xymatrix{
{ \operatorname{TR}_{q-\lambda}^r(A) } \ar[r]^(.42){\gamma}
\ar[d]^{R_s} & 
{ \prod_j \operatorname{TR}_{q-\lambda'}^u(A;p) }
\ar[d]^{R_s^{\gamma}} && 
{ \operatorname{TR}_{q-\lambda}^r(A) } \ar[r]^(.42){\gamma}
\ar[d]<1ex>^{F_s} & 
{ \prod_j \operatorname{TR}_{q-\lambda'}^u(A;p) }
\ar[d]<1ex>^{F_s^{\gamma}} \cr 
{ \operatorname{TR}_{q-\lambda'}^t(A) } \ar[r]^(.42){\gamma} & 
{ \prod_j \operatorname{TR}_{q-\lambda''}^{u-v}(A;p), } &&
{ \operatorname{TR}_{q-\lambda}^t(A) } \ar[r]^(.42){\gamma}
\ar[u]<1ex>^{V_s} & 
{ \prod_j \operatorname{TR}_{q-\lambda'}^{u-v}(A;p) }
\ar[u]<1ex>^{V_s^{\gamma}} \cr 
}$$
where the maps $R_s^{\gamma}$, $F_s^{\gamma}$, and $V_s^{\gamma}$ are
defined as follows: The map $R_s^{\gamma}$ takes the factor indexed by
a divisor $j$ of $t'$ to the factor indexed by the same divisor $j$ of
$t'$ by the map $R^v$ and annihilates the factors indexed by divisors
$j$ of $r'$ that do not divide $t'$. The map $F_s^{\gamma}$ takes the
factor indexed by a divisor $j$ of $r'$ that is divisible by $s'$ to
the factor indexed by the divisor $j/s'$ of $t'$ by the map $F^v$ and
annihilates the remaining factors. Finally, the map $V_s^{\gamma}$
takes the factor indexed by the divisor $j$ of $t'$ to the factor
indexed by the divisor $s'j$ of $r'$ by the map $s'V^v$.

It is now straightforward to check that the $p$-typical decomposition
of the top long-exact sequence in the diagram in the statement of
Thm.~\ref{maintr} takes the form
\begin{equation}\label{m}
\begin{aligned}
\cdots & \to
\prod_{j \in m'I_p} \lim_R \operatorname{TR}_{q-\lambda_d}^{u-v}(A;p)
\xrightarrow{m'V^v} 
\prod_{j \in I_p} \lim_R \operatorname{TR}_{q-\lambda_d}^u(A;p) \cr
{} & \xrightarrow{\varepsilon} \;\;\; K_{q+1}(A[x]/(x^m),(x)) \;\;\;
\xrightarrow{\partial} \,\, 
\prod_{j \in m'I_p} \lim_R
\operatorname{TR}_{q-1-\lambda_d}^{u-v}(A;p) \to \cdots \cr 
\end{aligned}
\end{equation}
where $m = p^v m'$ with $m'$ not divisible by $p$, and where $d =
d_p(m,u,j)$ is defined by
\begin{equation}\label{d_p(m,u,j)}
d_p(m,u,j) = \left[ \frac{p^{u-1}j - 1}{m} \right].
\end{equation}
In the top line, the right-hand product ranges over the set $I_p$ of
positive integers not divisible by $p$, and the left-hand product
ranges over the subset $m'I_p \subset I_p$. Similarly, the bottom
long-exact sequence in the diagram in the statement of
Thm.~\ref{maintr} takes the form
\begin{equation}\label{n}
\begin{aligned}
\cdots & \to
\prod_{j \in n'I_p} \lim_R \operatorname{TR}_{q-\lambda_e}^{u-w}(A;p)
\xrightarrow{n'V^w} 
\prod_{j \in I_p} \lim_R \operatorname{TR}_{q-\lambda_e}^u(A;p) \cr
{} & \xrightarrow{\varepsilon} \;\;\; K_{q+1}(A[x]/(x^n),(x)) \;\;\;
\xrightarrow{\partial} 
\, \, \prod_{j \in n'I_p} \lim_R
\operatorname{TR}_{q-1-\lambda_e}^{u-w}(A;p) \to \cdots \cr 
\end{aligned}
\end{equation}
where $n = p^wn'$ with $n'$ not divisible by $p$, and where $e =
d_p(n,u,j)$. It follows that Thm.~\ref{maintr} above is
equivalent to the following statement:

\begin{theorem}\label{ptypicalmaintr}Let $A$ be an
$\mathbb{F}_p$-algebra, and let $f \colon A[x]/(x^m) \to A[x]/(x^n)$
be the canonical projection. Then there is a map of long-exact
sequence from the sequence~(\ref{m}) to the sequence~(\ref{n}) that is
given, on the lower left-hand terms, by the map
$$f_* \colon K_*(A[x]/(x^m),(x)) \to K_*(A[x]/(x^n),(x))$$
induced by the canonical projection, on the upper right-hand terms, by
the map that takes the $j$th factor of the domain to the $j$th factor
of the target by the map
$$\iota_p(m,n,j)_q \colon \lim_R\operatorname{TR}_{q-\lambda_d}^u(A;p)
\to \lim_R\operatorname{TR}_{q-\lambda_e}^u(A;p)$$
induced from the canonical inclusions $\iota_p(m,n,j,u) \colon
S^{\lambda_d} \to S^{\lambda_e}$, and, on the upper left-hand terms,
by the zero map. The corresponding statement for the homotopy groups
with $\mathbb{Z}/p^i\mathbb{Z}$-coefficients is valid for every ring
$A$.
\end{theorem}

The following result shows that the limits that appear in the
statements of Thms.~\ref{maintr} and~\ref{ptypicalmaintr}
stabilize. In particular, the corresponding derived limits vanish. It
also shows that, for every integer $q$, the products that appear
in~(\ref{m}) and~(\ref{n}) are finite.

\begin{lemma}\label{limitreached}Let $u,u' \geqslant 1$ and $m > n
\geqslant 1$ be integers, and let $j$ be an integer not divisible by
$p$. Let $d = d_p(m,u,j)$ and $d' = d_p(m,u',j)$, and let $e =
d_p(n,u,j)$ and $e' = d_p(n,u',j)$. Then, in the commutative diagram,
$$\xymatrix{
{ \lim_R^{\phantom{R}}\operatorname{TR}_{q - \lambda_d}^u(A;p) }
\ar[rr]^{\iota_p(m,n,j)_q} \ar[d]^{\operatorname{pr}_{u'}} &&
{ \lim_R^{\phantom{R}}\operatorname{TR}_{q - \lambda_e}^u(A;p) }
\ar[d]^{\operatorname{pr}_{u'}} \cr 
{ \operatorname{TR}_{q - \lambda_{d'}}^{u'}(A;p) }
\ar[rr]^{\iota_p(m,n,j,u')_q} && 
{ \operatorname{TR}_{q - \lambda_{e'}}^{u'}(k;p), } \cr
}$$
the left-hand vertical map is an isomorphism, if $q < 2d_p(m,u'+1,j)$,
and right-hand vertical map is an isomorphism, if $q <
2d_p(n,u'+1,j)$.
\end{lemma}

\begin{proof}Let $r$ be a positive integer, let $\tilde{d} =
d_p(m,r,j)$, and let $\tilde{d}' = d_p(m,r-1,j)$. It follows
from~\cite[Thm.~2.2]{hm} that there is a long-exact sequence
$$\cdots \to
\mathbb{H}_q(C_{p^{r-1}},T(A) \wedge S^{\lambda_{\tilde{d}}})
\xrightarrow{N} 
\operatorname{TR}_{q-\lambda_{\tilde{d}}}^r(A;p) \xrightarrow{R} 
\operatorname{TR}_{q - \lambda_{\tilde{d}'}}^{r-1}(A;p) \to \cdots.$$
The Borel homology group that appear on the left-hand side is the
abutment of a first quadrant homology type spectral sequence
$$E_{s,t}^2 =
H_s(C_{p^{r-1}},\operatorname{TR}_{t - \lambda_{\tilde{d}}}^1(A;p)) 
\Rightarrow 
\mathbb{H}_{s+t}(C_{p^{r-1}}, T(A) \wedge S^{\lambda_{\tilde{d}}}),$$
and the groups $\operatorname{TR}_{q - \lambda_{\tilde{d}}}^1(A;p)$
are zero, for $q < 2\tilde{d}$. Hence, the map $R$ in the long-exact
sequence above is an isomorphism, for $q < 2\tilde{d}$. The lemma
follows.
\end{proof}

\section{The cyclic bar-construction}\label{barconstruction}

In this section, we prove Thm.~\ref{maintr} above. We first recall the
structure of the topological Hochschild $\mathbb{T}$-spectrum of the
truncated polynomial algebra $A[x]/(x^m)$. The reader is referred
to~\cite{h3} for details on the topological Hochschild
$\mathbb{T}$-spectrum and the cyclic bar-construction. We will write
$X[-]$ to indicate a cyclic object with $r$-simplices $X[r]$.

We showed in~\cite[Thm.~7.1]{hm} (see also~\cite[Prop.~4]{h3}) that
there is an $\mathscr{F}$-equivalence
$$\alpha \colon N^{\operatorname{cy}}(\Pi_m) \wedge T(A) \to
T(A[x]/(x^m))$$ 
where $N^{\operatorname{cy}}(\Pi_m)$ is the geometric realization of
the cyclic bar-construction $N^{\operatorname{cy}}(\Pi_m)[-]$ of the
pointed monoid $\Pi_m = \{0,1,x,x^2,\dots,x^{m-1}\}$ with base point
$0$ and with $x^m = 0$. A map of $\mathbb{T}$-spectra is an 
$\mathscr{F}$-equivalence, if it induces an equivalence of $C_r$-fixed 
point spectra, for all finite subgroups $C_r \subset \mathbb{T}$, and 
hence, we have induced isomorphisms
$$[ S^q \wedge (\mathbb{T}/C_r)_+,
N^{\operatorname{cy}}(\Pi_m) \wedge T(A) ]_{\mathbb{T}}
\xrightarrow{\alpha_*} 
[ S^q \wedge (\mathbb{T}/C_r)_+, T(A[x]/(x^m)) ]_{\mathbb{T}} 
 = \operatorname{TR}_q^r(A[x]/(x^m)).$$
Moreover, the following diagram commutes
$$\xymatrix{
{ N^{\operatorname{cy}}(\Pi_m) \wedge T(A) } \ar[r]^{\alpha}
\ar[d]^{f_*' \wedge \operatorname{id}} & 
{ T(A[x]/(x^m)) } \ar[d]^{f_*} \cr
{ N^{\operatorname{cy}}(\Pi_n) \wedge T(A) } \ar[r]^{\alpha} &
{ T(A[x]/(x^m)), } \cr
}$$
where $f \colon A[x]/(x^m) \to A[x]/(x^n)$ and $f' \colon \Pi_m \to
\Pi_n$ are the canonical projections. 

The cyclic bar-construction of $\Pi_m$ decomposes as the wedge sum of
pointed cyclic sets
$$N^{\operatorname{cy}}(\Pi_m)[-] =
\bigvee_{i \geqslant 0} N^{\operatorname{cy}}(\Pi_m,i)[-],$$
where $N^{\operatorname{cy}}(\Pi_m,i)[-] \subset
N^{\operatorname{cy}}(\Pi_m)[-]$ is the sub-pointed cyclic set whose
$k$-simplices are the $(k+1)$-tuples $(x^{i_0},\dots,x^{i_k})$, where
$i_0 + \dots + i_k = i$, and the base-point $0$. The geometric
realization decomposes accordingly as a wedge sum of pointed
$\mathbb{T}$-spaces
$$N^{\operatorname{cy}}(\Pi_m) = 
\bigvee_{i \geqslant 0} N^{\operatorname{cy}}(\Pi_m,i).$$ 
The homotopy type of the pointed $\mathbb{T}$-space
$N^{\operatorname{cy}}(\Pi_m,i)$ was  determined
in~\cite[Thm.~B]{hm1}. The result is that
$N^{\operatorname{cy}}(\Pi_m,0)$ is the discrete space $\{0,1\}$ and
that for positive $i$ there is a canonical cofibration sequence of
pointed
$\mathbb{T}$-spaces
$$\mathbb{T}_+ \wedge_{C_{i/m}} S^{\lambda_d}
\xrightarrow{\operatorname{pr}} 
\mathbb{T}_+ \wedge_{C_i} S^{\lambda_d} \xrightarrow{\varepsilon}
N^{\operatorname{cy}}(\Pi_m,i) \xrightarrow {\partial}
\Sigma \mathbb{T}_+ \wedge_{C_{i/m}} S^{\lambda_d},$$
where $d = d(m,i) = [(i-1)/m]$. The left-hand term is understood to be
a point, if $m$ does not divide $i$. More precisely, there is a
canonical homotopy class of maps of pointed $\mathbb{T}$-spaces
$$\hat{\theta}_d \colon N^{\operatorname{cy}}(\Pi_m,i) \to
\operatorname{cone}(
\mathbb{T}_+ \wedge_{C_{i/m}} S^{\lambda_d}
\xrightarrow{\operatorname{pr}} 
\mathbb{T}_+ \wedge_{C_i} S^{\lambda_d})$$
from $N^{\operatorname{cy}}(\Pi_m,i)$ to the mapping cone of the map
$\operatorname{pr}$ and any map in this homotopy class is a weak
equivalence of pointed $\mathbb{T}$-spaces. The canonical projection
$f' \colon \Pi_m \to \Pi_n$ induces a map of cofibration sequences
which we identify in Prop.~\ref{trianglemap} below. But first, we
recall from~\cite[\S3]{hm1} how the homotopy class of maps
$\hat{\theta}_d$ is defined.

We view the standard simplex $\Delta^{i-1}$ as the convex hull of the
set of group elements $C_i$ inside the regular representation
$\mathbb{R}[C_i]$. Let $\zeta_i \in C_i \subset \mathbb{T}$ be the
generator $\zeta_i = \exp(2 \pi \sqrt{-1}/i)$, and let $\Delta^{i-m}
\subset \Delta^{i-1}$ be the convex hull of the group elements 
$1, \zeta_i, \dots, \zeta_i^{i-m}$. As a cyclic set,
$N^{\operatorname{cy}}(\Pi_m,i)[-]$ is generated by the single
$(m-1)$-simplex $(x,\dots,x)$. Based on this observation, we showed
in~\cite[Lemma~2.2.6]{hm1} that there is a $\mathbb{T}$-equivariant
homeomorphism
$$\varphi \colon \mathbb{T}_+ \wedge_{C_i} (\Delta^{i-1}/C_i \cdot \Delta^{i-m})
\xrightarrow{\sim} N^{\operatorname{cy}}(\Pi_m,i).$$
The regular representation $\mathbb{R}[C_i]$ has the canonical direct sum
decomposition
$$\mathbb{R}[C_i] = \begin{cases}
\mathbb{C}(s) \oplus \dots \oplus \mathbb{C}(1) \oplus \mathbb{R}, &
\text{if $i = 2s+1$,} \cr
\mathbb{C}(s) \oplus \dots \oplus \mathbb{C}(1) \oplus \mathbb{R}
\oplus \mathbb{R}_{-}, & \text{if $i = 2s+2$,} \cr 
\end{cases}$$
so if $d \leqslant s$, or equivalently, if $2d < i$, there is a
canonical projection
$$\pi_d \colon \mathbb{R}[C_i] \to \mathbb{C}(d) \oplus \dots \oplus
\mathbb{C}(1) = \lambda_d.$$
Suppose first that $md < i < m(d+1)$. We showed
in~\cite[Thm.~3.1.2]{hm1} that, in this case, the image of $C_i \cdot
\Delta^{i-m} \subset \mathbb{R}[C_i]$ does not contain the origin 
$0 \in \lambda_d$. Hence, we may compose $\pi_d$ with the radial
projection away from a small ball around $0 \in \lambda_d$ to get a
$C_i$-equivariant map
$$\theta_d \colon \Delta^{i-1}/C_i \cdot \Delta^{i-m} \to
D(\lambda_d) / S(\lambda_d) = S^{\lambda_d}$$
whose homotopy class is well-defined. The composite map of pointed
$\mathbb{T}$-spaces
$$\xymatrix{
{ \hat{\theta}_d \colon N^{\operatorname{cy}}(\Pi_m,i) }
\ar[r]^(.43){\varphi^{-1}} & 
{ \mathbb{T}_+ \wedge_{C_i} (\Delta^{i-1} / C_i \cdot \Delta^{i-m}) }
\ar[r]^(.63){\operatorname{id} \wedge \theta_d} &
{ \mathbb{T}_+ \wedge_{C_i} S^{\lambda_d} } \cr
}$$
is then the desired map $\hat{\theta}_d$. Its homotopy class is
well-defined and~\cite[Prop.~3.2.6]{hm1} shows that it is a weak
equivalence.

Suppose next that $i = m(d+1)$. We let $\lambda_d' \subset
\lambda_{d+1}$ be the image of the canonical inclusion $\iota \colon
\lambda_d \to \lambda_{d+1}$ and let $\lambda_d^{\perp} \subset
\lambda_{d+1}$ be the orthogonal complement of $\lambda_d' \subset
\lambda_{d+1}$. The canonical projection from $\lambda_{d+1}$ to
$\mathbb{C}(d+1)$ restricts to an isomorphism of $\lambda_d^{\perp}$
onto $\mathbb{C}(d+1)$, and we define $C_m' \subset \lambda_d^{\perp}$
to be the preimage by this isomorphism of $C_m \subset
\mathbb{C}(d+1)$. It follows again from~\cite[Thm.~3.1.2]{hm1} that
the image $\pi_{d+1}(C_i \cdot \Delta^{i-m})$ does not contain the
origin $0 \in \lambda_{d+1}$ and, in addition, \cite[3.3.5]{hm1} shows
that
$$\pi_{d+1}(C_i \cdot \Delta^{i-m}) \cap \lambda_d^{\perp} = C_m'.$$
We pick a small open ball $B \subset \lambda_{d+1} \smallsetminus
C_m'$ around a point in $S(\lambda_d^{\perp})$ and define
$$U = (C_i \cdot B) \cap S(\lambda_{d+1}) \subset S(\lambda_{d+1}).$$
If the ball $B$ is small enough, then the composition of the
projection $\pi_{d+1}$ and the radial projection away from a 
(different) small ball around $0 \in \lambda_{d+1}$ defines a
$C_i$-equivariant map
$$\theta_d' \colon \Delta^{i-1} / C_i \cdot \Delta^{i-m} \to
D(\lambda_{d+1}) / ( S(\lambda_{d+1}) \smallsetminus U )$$
whose homotopy class is well-defined. Let $C(C_m) = \{0\} * C_m
\subset D(\mathbb{C}(d+1))$ be the unreduced cone with base
$C_m \subset \mathbb{C}(d+1)$ and with cone point $0 \in
\mathbb{C}(d+1)$. Then the canonical homeomorphism of
$D(\mathbb{C}(d+1)) \times D(\lambda_d)$ onto $D(\lambda_{d+1})$
induces an inclusion into the target of the map $\theta_d'$ of the
pointed $C_i$-space
$$\frac{C(C_m) \times D(\lambda_d)}{C(C_m) \times S(\lambda_d) \cup
C_m \times D(\lambda_{d})} = (S^0 * C_m) \wedge S^{\lambda_d}.$$
One immediately verifies that this inclusion is a strong deformation
retract of pointed $C_i$-spaces. Hence, the map $\theta_d'$ defines a
homotopy class of maps of pointed $C_i$-spaces
$$\theta_d \colon \Delta^{i-1} / C_i \cdot \Delta^{i-m} \to
(S^0 * C_m) \wedge S^{\lambda_d}.$$
The composite map of pointed $\mathbb{T}$-spaces
$$\xymatrix{
{ \hat{\theta}_d \colon N^{\operatorname{cy}}(\Pi_m,i) }
\ar[r]^(.43){\varphi^{-1}} &
{ \mathbb{T}_+ \wedge_{C_i} ( \Delta^{i-1} / C_i \cdot \Delta^{i-m} )} 
\ar[r]^(.49){\operatorname{id} \wedge \theta_d} &
{ \mathbb{T}_+ \wedge_{C_i} ( (S^0 * C_m) \wedge S^{\lambda_d} ) } \cr
}$$
is the desired map $\hat{\theta}_d$. Its homotopy class is
well-defined and~\cite[Prop.~3.3.9]{hm1} shows that it is a weak
equivalence. In preparation for the proof of Prop.~\ref{trianglemap}
below, we first prove the following key result.

\begin{lemma}\label{homotopycommutative}Let $m$ and $n$ be positive
integers with $m > n$. 

(i) If $md < i < m(d+1)$ and $ne < i < n(e+1)$, then $d \leqslant e$
and the diagram
$$\xymatrix{
{ \Delta^{i-1} / C_i \cdot \Delta^{i-m} } \ar[r]^(.67){\theta_d}
\ar[d]^{\operatorname{pr}} & 
{ S^{\lambda_d} } \ar[d]^{\iota} \cr
{ \Delta^{i-1} / C_i \cdot \Delta^{i-n} } \ar[r]^(.67){\theta_e} &
{ S^{\lambda_e} } \cr
}$$
commutes up to $C_i$-equivariant homotopy. Here the right-hand
vertical map $\iota$ is the canonical inclusion of $S^{\lambda_d}$ in
$S^{\lambda_e}$.

(ii) If $md < i < m(d+1)$ and $i = n(e+1)$, then $d \leqslant e$ and
the diagram
$$\xymatrix{
{ \Delta^{i-1} / C_i \cdot \Delta^{i-m} } \ar[r]^(.65){\theta_d}
\ar[d]^{\operatorname{pr}} & 
{ S^{\lambda_d} } \ar[d]^{\iota} \cr
{ \Delta^{i-1} / C_i \cdot \Delta^{i-n} } \ar[r]^{\theta_e} &
{ (S^0 * C_n) \wedge S^{\lambda_e} } \cr
}$$
commutes up to $C_i$-equivariant homotopy. Here the right-hand
vertical map $\iota$ is the canonical inclusion of $S^{\lambda_d}$ in
$S^{\lambda_e}$ followed by the canonical inclusion of $S^{\lambda_e}
= S^0 \wedge S^{\lambda_e}$ in $(S^0 * C_n) \wedge S^{\lambda_e}$.

(iii) If $i = m(d+1)$ and $ne < i < n(e+1)$, then $d < e$ and
the diagram
$$\xymatrix{
{ \Delta^{i-1} / C_i \cdot \Delta^{i-m} } \ar[r]^{\theta_d}
\ar[d]^{\operatorname{pr}} & 
{ ( S^0 * C_m ) \wedge S^{\lambda_d} } \ar[d]^{\iota} \cr
{ \Delta^{i-1} / C_i \cdot \Delta^{i-n} } \ar[r]^(.67){\theta_e} &
{ S^{\lambda_e} } \cr
}$$
commutes up to $C_i$-equivariant homotopy. Here the right-hand
vertical map $\iota$ is the canonical inclusion of $(S^0 * C_m)
\wedge S^{\lambda_d}$ in $(S^0 * S(\mathbb{C}(d+1))) \wedge
S^{\lambda_d} = S^{\lambda_{d+1}}$ followed by the canonical inclusion
of $S^{\lambda_{d+1}}$ in $S^{\lambda_e}$. 

(iv) If $i = m(d+1) = n(e+1)$, then $d < e$ and the diagram
$$\xymatrix{
{ \Delta^{i-1} / C_i \cdot \Delta^{i-m} } \ar[r]^{\theta_d}
\ar[d]^{\operatorname{pr}} & 
{ (S^0 * C_m) \wedge S^{\lambda_d} } \ar[d]^{\iota} \cr
{ \Delta^{i-1} / C_i \cdot \Delta^{i-n} } \ar[r]^{\theta_e} &
{ (S^0 * C_n) \wedge S^{\lambda_e} } \cr
}$$
commutes up to $C_i$-equivariant homotopy. Here the right-hand
vertical map $\iota$ is the canonical inclusion of $(S^0 * C_m)
\wedge S^{\lambda_d}$ in $(S^0 * S(\mathbb{C}(d+1))) \wedge
S^{\lambda_d} = S^{\lambda_{d+1}}$ followed by the canonical
inclusions of $S^{\lambda_{d+1}}$ in $S^{\lambda_e}$ and of
$S^{\lambda_e} = S^0 \wedge S^{\lambda_e}$ in $(S^0 * C_n) \wedge
S^{\lambda_e}$.
\end{lemma}

\begin{proof}We first consider the case~(i). Let $\lambda_e -
\lambda_d \subset \lambda_e$ be the orthogonal complement of
$\lambda_d \subset \lambda_e$. Then there is a canonical
$C_i$-equivariant homeomorphism
$$D(\lambda_e - \lambda_d) \times D(\lambda_d) \xrightarrow{\sim}
D(\lambda_e)$$
which restricts to a $C_i$-equivariant homeomorphism
$$D(\lambda_e - \lambda_d) \times S(\lambda_d) \cup
S(\lambda_e - \lambda_d) \times D(\lambda_d) \xrightarrow{\sim}
S(\lambda_e).$$
We note that the composition
$$\Delta^{i-1} / C_i \cdot \Delta^{i-m}
\xrightarrow{\operatorname{pr}} 
\Delta^{i-1} / C_i \cdot \Delta^{i-n} \xrightarrow{\theta_e}
D(\lambda_e) / S(\lambda_e)$$
of the left-hand vertical map and the lower horizontal map in the
diagram of part~(i) of the statement factors through the map
$$\frac{D(\lambda_e - \lambda_d) \times D(\lambda_d)}{
D(\lambda_e - \lambda_d)\times S(\lambda_d)} \to
D(\lambda_e) / S(\lambda_e)$$
induced by the canonical homeomorphism above. But the latter map is
homotopic to the composite map
$$\frac{D(\lambda_e - \lambda_d) \times D(\lambda_d)}{ 
D(\lambda_e - \lambda_d)\times S(\lambda_d)}
\xrightarrow{\operatorname{pr}_2}
D(\lambda_d) / S(\lambda_d) \xrightarrow{\iota} 
D(\lambda_e) / S(\lambda_e)$$
by the $C_i$-equivariant homotopy induced from the radial contraction
$$h \colon D(\lambda_e - \lambda_d) \times [0,1] \to D(\lambda_e -
\lambda_d)$$
defined by $h(z,t) = tz$. This homotopy, in turn, induces a
$C_i$-equivariant homotopy from the composite $\iota \circ \theta_d$
to the composite $\theta_e \circ \operatorname{pr}$ in the diagram in
part~(i) of the statement. 

In the case~(ii), we note similarly that the composition
$$\begin{aligned}
{} & \Delta^{i-1}/C_i \cdot \Delta^{i-m}
\xrightarrow{\operatorname{pr}} 
\Delta^{i-1}/C_i \cdot \Delta^{i-n} 
{} & \xrightarrow{\theta_e} \frac{C(C_n) \times D(\lambda_e)}{ 
C(C_n) \times S(\lambda_e) \cup  C_n\times D(\lambda_e)} \cr
\end{aligned}$$
of the left-hand vertical map and the lower horizontal map in the
diagram in part~(ii) of the statement factors through the canonical
projection
$$\frac{C(C_n) \times D(\lambda_e)}{C(C_n) \times S(\lambda_e)} \to
\frac{C(C_n) \times D(\lambda_e)}{C(C_n) \times S(\lambda_e) \cup
C_n \times D(\lambda_e)}$$
Again, the latter map is homotopic to the composition
$$\frac{C(C_n) \times D(\lambda_e)}{C(C_n) \times S(\lambda_e)}
\xrightarrow{\operatorname{pr}_2} D(\lambda_e) / S(\lambda_e)
\xrightarrow{\iota} \frac{C(C_n) \times D(\lambda_e)}{C(C_n) \times
S(\lambda_e) \cup C_n \times D(\lambda_e)}$$
by the $C_i$-equivariant homotopy induced from the radial contraction
$$h \colon C(C_n) \times [0,1] \to C(C_n)$$
defined by $h(z,t) = tz$. This homotopy induces the desired homotopy
from the composite map $\iota \circ \theta_d$ to the composite map
$\theta_e \circ \operatorname{pr}$ in the diagram in part~(ii) of the
statement.

In the case~(iii) of the statement, one proves as in the proof of
part~(i) of the statement that the composite map $\operatorname{pr}
\circ \theta_e$ 
is homotopic to the composition
$$\Delta^{i-1}/C_i \cdot \Delta^{i-m} \xrightarrow{\theta_{d+1}}
D(\lambda_{d+1}) / S(\lambda_{d+1}) \xrightarrow{\iota}
D(\lambda_e) / S(\lambda_e).$$
But the map $\theta_{d+1}$ is homotopic to the composition
$$\Delta^{i-1}/C_i \cdot \Delta^{i-m} \xrightarrow{\theta_d}
\frac{C(C_m) \times D(\lambda_d)}{C(C_m) \times S(\lambda_d) \cup
C_m \times D(\lambda_d)} 
\xrightarrow{\iota} D(\lambda_{d+1}) / S(\lambda_{d+1}).$$
Indeed, this is a direct consequence of the construction of the map
$\lambda_d$. This proves the case~(iii) of the statement. Finally, the
proof of the case~(iv) is analogous to the proof of the case~(i). This
completes the proof.
\end{proof}

\begin{prop}\label{trianglemap}Let $m > n \geqslant 1$ be integers and
$f' \colon \Pi_m \to \Pi_n$ the canonical projection. Let $i
\geqslant 1$ be an integer, and set $d = [(i-1)/m]$ and $e =
[(i-1)/n]$. Then there is a homotopy commutative diagram of pointed
$\mathbb{T}$-spaces
$$\xymatrix{
{ \mathbb{T}_+ \wedge_{C_{i/m}} S^{\lambda_d} }
\ar[r]^(.52){\operatorname{pr}} \ar[d]^{*} &
{ \mathbb{T}_+ \wedge_{C_i} S^{\lambda_d} }
\ar[r]^{\varepsilon} \ar[d]^{\operatorname{id} \wedge \iota} &
{ N^{\operatorname{cy}}(\Pi_m,i) } \ar[r]^(.42){\partial}
\ar[d]^{f_*'} &
{ \Sigma \mathbb{T}_+ \wedge_{C_{i/m}} S^{\lambda_d} }
\ar[d]^{*} \cr
{ \mathbb{T}_+ \wedge_{C_{i/n}} S^{\lambda_e} }
\ar[r]^(.52){\operatorname{pr}} &
{ \mathbb{T}_+ \wedge_{C_i} S^{\lambda_e} }
\ar[r]^{\varepsilon} &
{ N^{\operatorname{cy}}(\Pi_n,i) }
\ar[r]^(.42){\partial} &
{ \Sigma \mathbb{T}_+ \wedge_{C_{i/n}} S^{\lambda_e}, } \cr
}$$
where the rows are cofibration sequences, the map $\operatorname{id} \wedge \iota$
is induced from the canonical inclusion of $\lambda_d$ in $\lambda_e$,
and the map $*$ is the null-map. The domain (resp.~target) of the map
$*$ is understood to be a point if $m$ (resp.~$n$) does not divide
$i$.
\end{prop}

\begin{proof}As we recalled above, the statement that the rows in the
diagram of the statement are cofibration sequences of pointed
$\mathbb{T}$-spaces is equivalent to the statement that the maps of
$\mathbb{T}$-spaces $\hat{\theta}_d$ are weak equivalences, and the latter
statement is~\cite[Props.~3.2.6, 3.3.9]{hm1}. The statement
that the diagram commutes follows immediately from
Lemma~\ref{homotopycommutative} upon applying the functor that to a
pointed $C_i$-space $X$ associates the pointed $\mathbb{T}$-space
$\mathbb{T}_+ \wedge_{C_i} X$ with the exception of the statement that
the left-hand vertical map is null-homotopic. Only the case~(iii)
needs proof. In this case, we have a homotopy commutative diagram of
pointed $C_i$-spaces
$$\xymatrix{
{ C_{m+} \wedge S^{\lambda_d} } \ar[r]^(.6){\operatorname{pr}_2}
\ar[d] &
{ S^{\lambda_d} } \ar[r] \ar@{=}[d] &
{ (S^0 * C_m) \wedge S^{\lambda_d} } \ar[r]^{\partial} \ar[d] &
{ \Sigma C_{m+} \wedge S^{\lambda_d} } \ar[d] \cr
{ S(\mathbb{C}(d+1))_+ \wedge S^{\lambda_d} } \ar[r] &
{ S^{\lambda_d} } \ar[r]^(.47){\iota} &
{ S^{\lambda_{d+1}} } \ar[r]^(.36){\partial} &
{ \Sigma S(\mathbb{C}(d+1))_+ \wedge S^{\lambda_d} } \cr
}$$
which shows that the composite map $\iota \circ \operatorname{pr}_2$
is null-homotopic. It follows that the composite map of pointed
$\mathbb{T}$-spaces
$$\mathbb{T}_+ \wedge_{C_{i/m}} S^{\lambda_d}
\xrightarrow{\operatorname{pr}} 
\mathbb{T}_+ \wedge_{C_i} S^{\lambda_d}
\xrightarrow{\operatorname{id} \wedge \iota} 
\mathbb{T}_+ \wedge_{C_i} S^{\lambda_{d+1}}$$
is null-homotopic as desired.
\end{proof}

\begin{proof}[of Thm.~\ref{maintr}]The proof that the horizontal
sequences in the diagram in the statement are exact is given
in~\cite[Prop.~4.2.3]{hm}. A direct, and perhaps more detailed, proof
that the isomorphic sequences~(\ref{m}) and~(\ref{n}) are exact is
given in~\cite[Prop.~11]{h3}. It is immediately clear from either
proof that Prop.~\ref{trianglemap} implies that the diagram in the
statement commutes.
\end{proof}

\section{The map $\iota_* \colon
  \operatorname{TR}_{q-\lambda_{a-1}}^u(A;p) \to 
  \operatorname{TR}_{q-\lambda_a}^u(A;p)$}\label{FdV}

The topological Hochschild $\mathbb{T}$-spectrum $T(A)$ gives rise to
a generalized equivariant homology theory on the category of pointed
$\mathbb{T}$-spaces defined by
$$T(A)_q^{C_r}(X) = \pi_q^{C_r}(X \wedge T(A)) = 
[S^q \wedge (\mathbb{T}/C_r)_+, X \wedge
T(A)]_{\mathbb{T}},$$
and the map of the title is equal to the map of these homology groups
induced by the canonical inclusion $\iota \colon S^{\lambda_{a-1}} \to
S^{\lambda_a}$. We will evaluate this map in Prop.~\ref{exactsequence}
below but, in preparation, we first prove a general result about
$\mathbb{T}$-spectra. We define a $\mathbb{T}$-spectrum to be an
orthogonal $\mathbb{T}$-spectrum in the sense
of~\cite{mandellmay}. The category of $\mathbb{T}$-spectra is a model
category and the associated homotopy category, by definition, is the
$\mathbb{T}$-stable homotopy category. Moreover, the category of
$\mathbb{T}$-spectra has a closed symmetric monoidal structure given
by the smash product which induces a closed symmetric monoidal
structure on the $\mathbb{T}$-stable homotopy category. We give the 
$\mathbb{T}$-stable homotopy category the triangulated structure
defined in~\cite[\S2]{hm4}.

Let $T$ be a $\mathbb{T}$-spectrum, and let $C_r \subset \mathbb{T}$
be the subgroup of order $r$. We define
$$\pi_q^{C_r}(T) = [ S^q \wedge (\mathbb{T}/C_r)_+, T]_{\mathbb{T}}$$
to be the set of maps in the $\mathbb{T}$-stable homotopy category
between the indicated $\mathbb{T}$-spectra. In particular,
$\operatorname{TR}_{q-\lambda}^r(A) = \pi_q^{C_r}(S^{\lambda} \wedge
T(A))$. There are canonical isomorphisms
$$\pi_q(T^{C_r}) = 
[ S^q, T^{C_r}] \xrightarrow{\sim}
[ S^q \wedge (\mathbb{T}/C_r)_+, T^{C_r}]_{\mathbb{T}/C_r}
\xrightarrow{\sim} 
[S^q \wedge (\mathbb{T}/C_r)_+, T]_{\mathbb{T}} = \pi_q^{C_r}(T),$$
where, in the middle term, $T^{C_r}$ denotes the $C_r$-fixed point
$\mathbb{T}/C_r$-spectrum, and where, on the left-hand side, $T^{C_r}$ 
denotes the underlying non-equivariant spectrum of this
$\mathbb{T}/C_r$-spectrum. We consider the cofibration sequence of
pointed $\mathbb{T}$-spaces
$$\mathbb{T}_+ \xrightarrow{\pi} S^0 \xrightarrow{\iota}
S^{\mathbb{C}(1)} \xrightarrow{\partial} 
S^1 \wedge \mathbb{T}_+ = \Sigma \mathbb{T}_+,$$
where $\pi$ collapses $\mathbb{T}$ onto the non-base point of $S^0$,
and where $\iota$ is the canonical inclusion. We identify the
underlying pointed spaces of $S^2$ and $S^{\mathbb{C}(1)}$ by the
isomorphism
$$\varphi \colon S^2 \to S^{\mathbb{C}(1)}$$
that takes the class of $(a,b)$ in $S^2 = S^1 \wedge S^1$ to the class 
of $a + b\sqrt{-1}$ in $S^{\mathbb{C}(1)}$. The composition of
$\varphi$ and $\partial$ defines a map in the non-equivariant stable
homotopy category
$$\sigma \colon S^1 \to \mathbb{T}_+.$$
The cofibration sequence above induces a cofibration sequence in the
non-equivariant stable homotopy category, and the latter sequence
splits. We thus obtain a direct sum diagram in the non-equivariant
stable homotopy category
$$\xymatrix{
{ S_{\phantom{+}}^1 } \ar[r]<.6ex>^{\sigma} &
{ \mathbb{T}_+^{\phantom{1}} } \ar[r]<.6ex>^{\pi} \ar[l]<.6ex>^{\kappa} &
{ S_{\phantom{+}}^0 } \ar[l]<.6ex>^{i} \cr
}$$
where $i$ is the section of $\pi$ that takes the the non-base point in
$S^0$ to $1 \in \mathbb{T}$, and where the map $\kappa$ is the
corresponding retraction of the map $\sigma$. We define the map
$$d \colon \pi_q^{C_r}(T) \to \pi_{q+1}^{C_r}(T) \hskip7mm
\text{(Connes' operator)}$$
to be composition
$$\xymatrix{
{ \pi_q(T^{C_r}) } \ar[r]^(.4){\ell_{\sigma}} &
{ \pi_{q+1}(\mathbb{T}_+ \wedge T^{C_r}) } \ar[r]^(.46){\rho_{r*}} &
{ \pi_{q+1}((\mathbb{T}/C_r)_+ \wedge T^{C_r}) } \ar[r]^(.62){\mu_*} &
{ \pi_{q+1}(T^{C_r}), } \cr
}$$
where the left-hand map is left multiplication by $\sigma$, where the
middle map is induced by the isomorphism $\rho_r \colon \mathbb{T} \to
\mathbb{T}/C_r$ defined by the $r$th root, and where the right-hand
map is induced from the left action by $\mathbb{T}$ on the underlying
non-equivariant spectrum of $T$. It anti-commutes with the suspension
isomorphism in the sense that the following diagram anti-commutes:
$$\xymatrix{
{ \pi_q^{C_r}(T) } \ar[r]^(.48){d} \ar[d]^{\operatorname{susp}}
\ar@{}[dr]|-{(-1)} &
{ \pi_{q+1}^{C_r}(T) } \ar[d]^{\operatorname{susp}} \cr
{ \pi_{q+1}^{C_r}(\Sigma T) } \ar[r]^(.49){d} &
{ \pi_{q+2}^{C_r}(\Sigma T). } \cr
}$$
Moreover, $dd(x) = \eta \cdot d(x)$ and $F_rdV_r(x) =
d(x) + (r-1)\eta \cdot x$, where $\eta \in \pi_1(S^0)$ is the Hopf
class, and, if $T$ is a ring $\mathbb{T}$-spectrum, $d$ is a
derivation for the multiplication;~see~\cite[\S1]{h}. If $A$ is an
$\mathbb{F}_p$-algebra, then multiplication by $\eta$ is zero on
$\pi_*^{C_r}(T(A))$. Hence, in this case, Connes' operator is a
differential and satisfies $F_rdV_r = d$.

\begin{lemma}\label{FdVlemma}Let $T$ be a $\mathbb{T}$-spectrum, and
let $r$ and $j$ be relative prime positive integers. Then, for every
integer $q$, there is a commutative diagram
$$\xymatrix{
{ \pi_{q-1}(T) \oplus \pi_q(T) } \ar[r]^(.45){\operatorname{susp}
  \oplus \operatorname{susp}} \ar[d]^{i_* \oplus i_*} & 
{ \pi_q( S^1 \wedge T ) \oplus \pi_{q+1}( S^1 \wedge T ) } \cr
{ \pi_{q-1}((\mathbb{T}/C_j)_+ \wedge T) \oplus
\pi_q((\mathbb{T}/C_j)_+ \wedge T) } \ar[d]^{dV_r+V_r} &
{ \pi_q((\mathbb{T}/C_j)_+ \wedge T) \oplus
\pi_{q+1}((\mathbb{T}/C_j)_+ \wedge T) } \ar[u]_{\kappa_* \oplus
\kappa_*} \cr
{ \pi_q^{C_r}((\mathbb{T}/C_j)_+ \wedge T) } \ar@{=}[r] &
{ \pi_q^{C_r}((\mathbb{T}/C_j)_+ \wedge T ), } \ar[u]_{(F_r,F_rd)} \cr
}$$
where the top horizontal map is the suspension isomorphism. Moreover,
the compositions of the left-hand vertical maps and the right-hand
vertical maps are both isomorphisms.
\end{lemma}

\begin{proof}The general formulas that we recalled before the
statement show that
$$\big((\kappa_* \oplus \kappa_*) \circ (F_r,F_rd) \circ (dV_r+V_r)
\circ (i_* \oplus i_*) \big) (x,y) = \big( (\kappa_* \circ d \circ
i_*)(x), (\kappa_* \circ d \circ i_*)(y) \big),$$
and one verifies that $\kappa_* \circ d \circ i_*$ is equal to the
suspension isomorphism. This shows that the diagram commutes. It was
proved in~\cite[Prop.~3.4.1]{hm3} that the composition of the
left-hand vertical maps is an isomorphism, provided that $r = p^{v-1}$
is a power of a prime number $p$. The proof in the general case is
completely similar.
\end{proof}

Suppose now that $r$ is a divisor in $a$. Then the isomorphism
$$\varphi \colon S^2 \to S^{\mathbb{C}(a)}$$
is an isomorphism of pointed $C_r$-spaces. We define
$$\hat{\varphi} \colon S^2 \wedge S^q \wedge (\mathbb{T}/C_r)_+ \to
S^{\mathbb{C}(a)} \wedge S^q \wedge (\mathbb{T}/C_r)_+$$
be the isomorphism of pointed $\mathbb{T}$-spaces that takes the
class of $(x,y,zC_r)$ to the class of $(z\varphi(x),y,zC_r)$. It
follows that, if $r$ is a divisor in $a$, we have an isomorphism
$$\varphi^{!} \colon \operatorname{TR}_{q - \lambda_{a-1}}^r(A) \to
\operatorname{TR}_{q + 2 - \lambda_a}^r(A)$$
defined by the composition
$$\begin{aligned}
{} \operatorname{TR}_{q-\lambda_{a-1}}^r(A) & \;\; = \hskip7pt
[ S^q \wedge (\mathbb{T}/C_r)_+, S^{\lambda_{a-1}} \wedge T(A)
]_{\mathbb{T}} \cr
{} & \xrightarrow{\operatorname{susp}}
[ S^{\mathbb{C}(a)} \wedge S^q \wedge (\mathbb{T}/C_r)_+,
S^{\mathbb{C}(a)} \wedge S^{\lambda_{a-1}} \wedge T(A) ]_{\mathbb{T}} \cr
{} & \; \xrightarrow{\hat{\varphi}^*} \hskip3pt
[ S^2 \wedge S^q \wedge (\mathbb{T}/C_r)_+,
S^{\mathbb{C}(a)} \wedge S^{\lambda_{a-1}} \wedge T(A)]_{\mathbb{T}}
\;\; = \;\; \operatorname{TR}_{q+2-\lambda_a}^r(A)\cr 
\end{aligned}$$
We write $\varphi_{!}$ for the inverse of $\varphi^{!}$.

\begin{prop}\label{exactsequence}Let $A$ be a
$\mathbb{Z}_{(p)}$-algebra, let $a, u \geqslant 1$ be integers, and
let $v = v(a,u)$ be the minimum of $u$ and $v_p(a) + 1$. Then there
are natural long-exact sequences
$$\xymatrix{
{ : } \ar[d] &
{ : } \ar[d] \cr
{ \operatorname{TR}_{q-1-\lambda_{a-1}}^v(A;p) \oplus
\operatorname{TR}_{q-\lambda_{a-1}}^v(A;p) } 
\ar[d]^{dV^{u-v} + V^{u-v}} & 
{ \operatorname{TR}_{q+1-\lambda_a}^u(A;p) \oplus
\operatorname{TR}_{q-\lambda_{a-1}}^u(A;p) }
\ar[d]^{\operatorname{pr}_2} \cr 
{ \operatorname{TR}_{q-\lambda_{a-1}}^u(A;p) } \ar[d]^{\iota_*} &
{ \operatorname{TR}_{q-\lambda_{a-1}}^u(A;p) } \ar[d]^{\iota_* = 0}
\cr 
{ \operatorname{TR}_{q-\lambda_a}^u(A;p) } 
\ar[d]^{ (\varphi_{!} F^{u-v}, - \varphi_{!} F^{u-v}d) } &
{ \operatorname{TR}_{q-\lambda_a}^u(A;p) }
\ar[d]^{\operatorname{in}_1} \cr 
{ \operatorname{TR}_{q-2-\lambda_{a-1}}^v(A;p) \oplus
\operatorname{TR}_{q-1-\lambda_{a-1}}^v(A;p) } \ar[d] &
{ \operatorname{TR}_{q-\lambda_a}^u(A;p) \oplus
\operatorname{TR}_{q-1-\lambda_{a-1}}^u(A;p) } \ar[d] \cr
{ : } &
{ : } \cr
}$$
where the left-hand sequence is valid, for $v < u$, and
where the right-hand sequence is valid, for $v = u$.
\end{prop}

\begin{proof}We consider the cofibration sequence of pointed
$\mathbb{T}$-spaces
$$S(\mathbb{C}(a))_+ \xrightarrow{\pi} S^0 \xrightarrow{\iota}
S^{\mathbb{C}(a)} \xrightarrow{\delta} 
S^1 \wedge S(\mathbb{C}(a))_+$$
and the induced cofibration sequence in the $\mathbb{T}$-stable
category
$$S(\mathbb{C}(a))_+ \wedge T \xrightarrow{\pi'}
T \xrightarrow{\iota'}
S^{\mathbb{C}(a)} \wedge T \xrightarrow{\delta'}
S^1 \wedge S(\mathbb{C}(a))_+ \wedge T$$
where $T = S^{\lambda_{a-1}} \wedge T(A)$. We also abbreviate $r =
p^{u-1}$. The associated exact sequence of equivariant homotopy groups 
takes the form
$$\pi_q^{C_r}(S(\mathbb{C}(a))_+ \wedge T) \xrightarrow{\pi_*'}
\pi_q^{C_r}(T) \xrightarrow{\iota_*'}
\pi_q^{C_r}(S^{\mathbb{C}(a)} \wedge T) \xrightarrow{\partial}
\pi_{q-1}^{C_r}(S(\mathbb{C}(a))_+ \wedge T),$$
where $\partial$ is the composition of $\partial_*'$ and the inverse
of the suspension isomorphism. We identify this sequence with the
long-exact sequences of the statement. The terms
$\pi_q^{C_r}(T)$ and $\pi_q^{C_r}(S^{\mathbb{C}(a)} \wedge T)$ are
equal to $\operatorname{TR}_{q-\lambda_{a-1}}^u(A;p)$ and
$\operatorname{TR}_{q-\lambda_a}^u(A;p)$, respectively, and the map
$\iota_*'$ is equal to the map $\iota_*$. It remains to identify the
remaining term and the two maps $\pi_*'$ and $\partial_*$.

Suppose first that $u = v$ such that $\mathbb{C}(a)$ is
$C_r$-trivial. Then the map $i \colon S^0 \to S(\mathbb{C}(a))_+$ that
takes the non-base point in $S^0$ to $1 \in S(\mathbb{C}(a))$ is
$C_r$-equivariant and defines a section of the map $\pi_*'$ in the
exact sequence of equivariant homotopy groups above. This completes
the proof in the case $u = v$. 

We next assume that $v < u$ and consider first the case $v = 1$ such 
that $a$ and $r$ are relatively prime. The map $i \colon S^0 \to
S(\mathbb{C}(a))_+$ is a non-equivariant section of the projection
$\pi$. The following diagram commutes
$$\xymatrix{
{ \pi_q^{C_r}(S(\mathbb{C}(a))_+ \wedge T) } \ar[r]^(.6){\pi_*'} &
{ \pi_q^{C_r}(T) } \cr
{ \pi_{q-1}(S(\mathbb{C}(a))_+ \wedge T) \oplus
  \pi_q(S(\mathbb{C}(a))_+ \wedge T) } 
\ar[u]_{dV_r+V_r} \ar[r]^(.67){\pi_*' \oplus \pi_*'} &
{ \pi_{q-1}(T) \oplus \pi_q(T) } \ar[u]_{dV_r+V_r} \cr
{ \pi_{q-1}(T) \oplus \pi_q(T) } \ar[u]_{i_*' \oplus i_*'} &
{} \cr
}$$
and the composite map $(\pi_*' \oplus \pi_*') \circ (i_* \oplus i_*)$
is the identity. Since $S(\mathbb{C}(a))_+$ is isomorphic to
$(\mathbb{T}/C_a)_+$ as a pointed $\mathbb{T}$-space,
Lemma~\ref{FdVlemma} shows that the composition of the left-hand
vertical maps is an isomorphism. We use this isomorphism to identify
the left-hand term of the exact sequence of equivariant homotopy
groups above. The composition of this isomorphism and the map $\pi_*'$
is equal to the map
$$dV_r + V_r \colon \operatorname{TR}_{q-1-\lambda_{a-1}}^1(A;p) \oplus
\operatorname{TR}_{q-\lambda_{a-1}}^1(A;p) \to
\operatorname{TR}_{q-\lambda_{a-1}}^u(A;p)$$ 
as stated. Similarly, the following diagram commutes
$$\xymatrix{
{ \pi_q^{C_r}(S^{\mathbb{C}(a)} \wedge T) } \ar[r]^(.48){\partial_*}
\ar[d]^{(F_r,-F_rd)} &
{ \pi_{q-1}^{C_r}(S(\mathbb{C}(a))_+ \wedge T) } \ar[d]^{(F_r,F_rd)} \cr
{ \pi_q(S^{\mathbb{C}(a)} \wedge T) \oplus \pi_{q+1}(S^{\mathbb{C}(a)}
  \wedge T) } \ar[r]^(.45){\partial_* \oplus \partial_*} &
{ \pi_{q-1}(S(\mathbb{C}(a))_+ \wedge T) \oplus
  \pi_q(S(\mathbb{C}(a))_+ \wedge T). }  
\ar[d]^{\operatorname{susp}^{-1} \circ \kappa_*' \oplus 
\operatorname{susp}^{-1} \circ \kappa_*'} \cr
{} &
{ \pi_{q-2}(T) \oplus \pi_{q-1}(T) } \cr
}$$
and the composition of the lower horizontal map and the lower
right-hand vertical map is equal to the map $\varphi_{!} \oplus
\varphi_{!}$. Finally, Lemma~\ref{FdVlemma} shows that the
composition of the right-hand vertical maps is equal to the inverse of
the map $dV_ri_*\oplus V_ri_*$. This completes the proof in the case
where $1 = v < u$.

Finally, we suppose that $1 < v < u$ and abbreviate $r = p^{u-1}$, $t
= p^{v-1}$, and $s = p^{u-v}$. We recall that the root isomorphism
$\rho_t \colon \mathbb{T} \to \mathbb{T}/C_t$ induces an equivalence
of categories $\rho_t^*$ from the $\mathbb{T}/C_t$-stable category to
the $\mathbb{T}$-stable category. There is a commutative diagram
$$\xymatrix{
{ \pi_q^{C_r}(S(\mathbb{C}(a))_+ \wedge T) } \ar[r]^(.6){\pi_*'} &
{ \pi_q^{C_r}(T) } \ar[r]^(.42){\iota_*'} &
{ \pi_q^{C_r}(S^{\mathbb{C}(a)} \wedge T) } \ar[r]^(.4){\partial} &
{ \pi_{q-1}^{C_r}(S(\mathbb{C}(a))_+ \wedge T) } \cr
{ \pi_q^{C_s}(S(\mathbb{C}(j))_+ \wedge T') } \ar[r]^(.6){\pi_*'}
\ar[u] & 
{ \pi_q^{C_s}(T') } \ar[r]^(.42){\iota_*'} \ar[u] &
{ \pi_q^{C_s}(S^{\mathbb{C}(j)} \wedge T') } \ar[r]^(.4){\partial}
\ar[u] & 
{ \pi_{q-1}^{C_r}(S(\mathbb{C}(j))_+ \wedge T'), } \ar[u] \cr
}$$
where $T = S^{\lambda_{a-1}} \wedge T(A)$ and $T' =
\rho_t^*(T^{C_t})$, and where $j = a/t$. The vertical maps are the
composition of the isomorphism $\rho_t^* \colon
\pi_q^{C_r/C_t}(X^{C_t}) \to \pi_q^{C_s}(\rho_t^*(X^{C_t}))$ and the
canonical isomorphism of $\pi_q^{C_r}(X)$ and
$\pi_q^{C_r/C_t}(X^{C_t})$. The case of the lower sequence was treated
above. This completes the proof.
\end{proof}

\begin{remark}Let $A$ be a regular noetherian ring that is also an
$\mathbb{F}_p$-algebra. There is a statement similar to
Lemma~\ref{exactsequence} above for the $p$-typical de~Rham-Witt
complex of $A$. Suppose that $v < u$. Then there is a long-exact
sequence sequence 
$$\xymatrix{
{ \cdots } \ar[r]^(.3){\beta} &
{ W_v\Omega_A^{q-1} \oplus W_v\Omega_A^q } \ar[r]^(.64){\alpha} &
{ W_u\Omega_A^q } \ar[r]^{p^{u-v}} &
{ W_u\Omega_A^q } \ar[r]^(.35){\beta} &
{ W_v\Omega_A^q \oplus W_v\Omega_A^{q+1} } \ar[r]^(.67){\alpha} &
{ \cdots, } \cr
}$$
where $\alpha = dV^{u-v}+V^{u-v}$ and $\beta=(F^{u-v},-F^{u-v}d)$.
This sequence is a more precise version of the
statement~\cite[III (3.3.3.3)]{illusieraynaud}. We outline the steps
in the proof. First, the basic case $A = \mathbb{F}_p$ is clear. Next,
one uses~\cite[Thm.~B]{hm3} to show that the exactness of the sequence
for $A$ implies the exactness of the sequence for $A[x]$. Third, one
uses~\cite[Prop.~6.2.3]{hm3} or~\cite[Prop.~I.1.14]{illusie} to
conclude that the sequence is exact whenever $A$ is a smooth
$\mathbb{F}_p$-algebra. Finally, one uses the theorem of
Popescu~\cite{popescu} that every regular noetherian
$\mathbb{F}_p$-algebra is a filtered colimit of smooth
$\mathbb{F}_p$-algebras and the fact that, being an initial object,
the de~Rham-Witt complex commutes with colimits.
\end{remark}

\section{Proof of Thm.~\ref{main}}\label{proofoftheorem}

In this section we evaluate the maps $\iota_p(m,n,j)_q$ that appear in 
the statement of Thm.~\ref{ptypicalmaintr}. This proves
Thm.~\ref{ptypicalmain} and hence the equivalent Thm.~\ref{main}. We
first recall the following result which was proved
in~\cite[Prop.~9.1]{hm}.

\begin{prop}\label{TRperfect}Let $a$ and $u$ be positive
integers. Then for every non-negative integer $i$, there exists an
isomorphism
$$\sigma_p(a,u,i) \colon W_r(\mathbb{F}_p) \xrightarrow{\sim}
\operatorname{TR}_{2i-\lambda_a}^u(\mathbb{F}_p;p),$$
where the length $r = r_p(a,u,i)$ is given by
$$r_p(a,u,i) =
\begin{cases}
u, & \text{if $a \leqslant i$,} \cr
u - s, & \text{if $[a/p^s] \leqslant i < [a/p^{s-1}]$ and
$1 \leqslant s < u$,} \cr
0, & \text{if $i < [a/p^{u-1}]$.} \cr
\end{cases}$$
The group $\operatorname{TR}_{q-\lambda_a}^u(\mathbb{F}_p;p)$ is zero,
if $q$ is a negative or odd integer. In addition, the isomorphisms
$\sigma_p(a,u,i)$ may be chosen in such a way that the square diagrams
$$\xymatrix{
{ W_r(\mathbb{F}_p) } \ar[rr]^(.45){\sigma_p(a,u,i)} \ar[d]<1ex>^{F}
&& 
{ \operatorname{TR}_{2i-\lambda_a}^u(\mathbb{F}_p;p) } \ar[d]<1ex>^{F}
\cr 
{ W_{r-1}(\mathbb{F}_p) } \ar[rr]^(.45){\sigma_p(a,u-1,i)}
\ar[u]<1ex>^{V} && 
{ \operatorname{TR}_{2i-\lambda_a}^{u-1}(\mathbb{F}_p;p), }
\ar[u]<1ex>^{V} \cr 
}$$
where $r = r_p(a,u,i)$, commute.
\end{prop}

We recall the integer functions $s_p(m,i,j)$ and $d_p(m,u,j)$ defined
in~(\ref{s_p(m,i,j)}) and~(\ref{d_p(m,u,j)}).

\begin{cor}\label{limitreachedcor}Let $m,u,u' \geqslant 1$ and $i
\geqslant 0$ be integers, let $j \geqslant 1$ be an integer that is
not divisible by $p$, and let $d = d_p(m,u,j)$ and $d' =
d_p(m,u',j)$. Then the canonical projection
$$\operatorname{pr}_{u'} \colon
\lim_R\operatorname{TR}_{2i-\lambda_d}^u(\mathbb{F}_p;p) \to 
\operatorname{TR}_{2i-\lambda_{d'}}^{u'}(\mathbb{F}_p;p)$$
is an isomorphism, if $p^{u'}j > m(i+1)$, and the common group is a
cyclic $W(\mathbb{F}_p)$-module of length $s_p(m,i,j)$. The group
$\lim_R\operatorname{TR}_{q-\lambda_d}^u(\mathbb{F}_p;p)$ is zero, if
$q$ is a negative or odd integer.
\end{cor}

\begin{proof}We proved in Lemma~\ref{limitreached} that the map of the 
statement is an isomorphism, provided that $i < d_p(m,u'+1,j) =
[(p^{u'}j-1)/m]$. This constraint on the integer $u'$ is equivalent to
the inequality $i + 1 \leqslant [(p^{u'}j - 1)/m]$, which is
equivalent to $i + 1 \leqslant (p^{u'}j-1)/m$, which, in turn, is
equivalent to the stated inequality $m(i+1) < p^{u'}j$. Finally, a
similar calculation based on Prop.~\ref{TRperfect} shows that the
length of the $W(\mathbb{F}_p)$-module
$\lim_R\operatorname{TR}_{2i-\lambda_d}^u(\mathbb{F}_p;p)$ is equal to
$s_p(m,i,j)$ as stated.
\end{proof}

\begin{lemma}\label{preliminaryformula}Let $m > n \geqslant 1$,
$u \geqslant 1$, and $i \geqslant 0$ be integers. Let
$1 \leqslant j \leqslant m(i+1)$ be an integer that is not divisible
by $p$, and let $d = d_p(m,u,j)$ and $e = d_p(n,u,j)$. Then the map
$$\iota_p(m,n,j,u)_{2i} \colon
\operatorname{TR}_{2i-\lambda_d}^u(\mathbb{F}_p;p) \to 
\operatorname{TR}_{2i-\lambda_e}^u(\mathbb{F}_p;p)$$
takes a generator of the domain to the product of a generator of the
target and an element $\alpha' = \alpha_p'(m,n,i,j,u)$ of
$W(\mathbb{F}_p)$ with $p$-adic valuation 
$$v_p(\alpha') = \sum_{d < a \leqslant e}
\operatorname{length}_{W(\mathbb{F}_p)}
\operatorname{TR}_{2i-\lambda_a}^{v(u,a)}(\mathbb{F}_p;p),$$ 
where $v(u,a) = \min\{u, v_p(a)+1\}$.
\end{lemma}

\begin{proof}The canonical inclusion $\iota_p(m,n,j,u) \colon
\lambda_d \to \lambda_e$ is equal to the composition of the canonical
inclusions $\iota \colon \lambda_{a-1} \to \lambda_a$ for $d < a
\leqslant e$. Since the groups
$\operatorname{TR}_{q-\lambda_a}^v(\mathbb{F}_p;p)$ are zero, 
for $q$ odd, Prop.~\ref{exactsequence} above identifies the cokernel
of the map
$$\iota_* \colon
\operatorname{TR}_{2i-\lambda_{a-1}}^u(\mathbb{F}_p;p) \to 
\operatorname{TR}_{2i-\lambda_a}^u(\mathbb{F}_p;p)$$
with the $W(\mathbb{F}_p)$-module
$\operatorname{TR}_{2i-\lambda_a}^{v(u,a)}(\mathbb{F}_p;p)$. Equivalently,
$\iota_*$ takes a generator of the domain to the product of a
generator of the target and an element $\alpha_a' \in W(\mathbb{F}_p)$
of $p$-adic valuation 
$$v_p(\alpha_a') = \operatorname{length}_{W(\mathbb{F}_p)}
\operatorname{TR}_{2i-\lambda_a}^{v(u,a)}(\mathbb{F}_p;p).$$
The lemma follows.
\end{proof}

We proceed to manipulate the sum that appears in
Lemma~\ref{preliminaryformula}. To this end, we fix the positive
integer $u$, and consider the bi-graded $\mathbb{F}_p$-vector space
$E(u)$ defined by the associated graded for the $p$-adic filtration of
the $W(\mathbb{F}_p)$-modules that appear in
Lemma~\ref{preliminaryformula}, 
$$E(u)_{i,a} = \bigoplus_{r \geqslant 0} \operatorname{gr}_p^r
\operatorname{TR}_{2i-\lambda_a}^{v(u,a)}(\mathbb{F}_p;p).$$
The following results identifies the structure of this bi-graded
$\mathbb{F}_p$-vector space.

\begin{lemma}\label{gradedmodule}The bi-graded $\mathbb{F}_p$-vector
space $E(u)$ is isomorphic to the bi-graded $\mathbb{F}_p$-vector
space defined by the sum of symmetric algebras
$$A(u) = \bigoplus_{1 \leqslant r \leqslant u}
S_{\mathbb{F}_p}\{x_r,\sigma_r\},$$ 
where $\deg x_r = (p^{r-1},1)$ and $\deg \sigma_r = (0,1)$.
\end{lemma}

\begin{proof}The statement is equivalent to the equality
$$\dim_{\mathbb{F}_p} A(u)_{a,i} =
\operatorname{length}_{W(\mathbb{F}_p)}
\operatorname{TR}_{2i-\lambda_a}^{v(u,a)}(\mathbb{F}_p;p),$$ 
which we verify by direct calculation. If $a \leqslant i$, then
$A(u)_{a,i}$ has basis
$$x_1^a \sigma_1^{i-a}, x_2^{a/p} \sigma_2^{i - a/p}, \dots,
x_{v(u,a)}^{a/p^{v(u,a)-1}} \sigma_{v(u,a)}^{i - a/p^{v(u,a)-1}}$$
so $\dim_{\mathbb{F}_p} A(u)_{a,i} = v(u,a)$ as required. Similarly,
if $a/p^s \leqslant i < a/p^{s-1}$ with $1 \leqslant s < v(u,a)$, then 
$A(u)_{a,i}$ has basis
$$x_{s+1}^{a/p^s} \sigma_{s+1}^{i - a/p^s}, x_{s+2}^{a/p^{s+1}}
\sigma_{s+2}^{i - a/p^{s+1}}, \dots,
x_{v(u,a)}^{a/p^{v(u,a)-1}} \sigma_{v(u,a)}^{i - a/p^{v(u,a)-1}}$$
which shows that $\dim_{\mathbb{F}_p} A(u)_{a,i} = v(u,a)-s$ as
desired. Finally, if $i < a/p^{v(u,a)-1}$, then $A(u)_{a,i} = 0$. This
completes the proof.
\end{proof}

\begin{prop}\label{finalformula}Let $m > n \geqslant 1$ and $i
\geqslant 0$ be integers, and let $1 \leqslant j \leqslant m(i+1)$ be
an integer that is not divisible by $p$. Then the map
$$\iota_p(m,n,j)_{2i} \colon
\lim_R\operatorname{TR}_{2i-\lambda_d}^u(\mathbb{F}_p;p) \to 
\lim_R\operatorname{TR}_{2i-\lambda_e}^u(\mathbb{F}_p;p),$$
where $d = d_p(m,u,j)$ and $e = d_p(n,u,j)$, takes a generator of the 
domain to the product of a generator of the target and an element
$\alpha = \alpha_p(m,n,i,j)$ of $W(\mathbb{F}_p)$ with $p$-adic
valuation  
$$v_p(\alpha) = \sum_{0 \leqslant h < i} \big( s_p(m,h,j) - s_p(n,h,j)
\big).$$
\end{prop}

\begin{proof}By Cor.~\ref{limitreachedcor} above, we may instead
consider the map
$$\iota_p(m,n,j,u)_{2i} \colon
\operatorname{TR}_{2i-\lambda_d}^u(\mathbb{F}_p;p) \to 
\operatorname{TR}_{2i-\lambda_e}^u(\mathbb{F}_p;p),$$
for a fixed positive integer $u$ with $p^uj > m(i+1)$. We showed in
Lemma~\ref{preliminaryformula} above that this map takes a generator
of the domain to the product of a generator of the target and an
element $\alpha' = \alpha_p'(m,n,i,j,u)$ of $W(\mathbb{F}_p)$ whose
$p$-adic valuation is given by the sum
$$v_p(\alpha') = \sum_{d < a \leqslant e}
\operatorname{length}_{W(\mathbb{F}_p)} 
\operatorname{TR}_{2i - \lambda_a}^{v(u,a)}(\mathbb{F}_p;p)$$
where $v(u,a) = \min\{u,v_p(a)+1\}$. We show that this sum is equal to
the sum $v_p(\alpha)$ of the statement. It follows from
Lemma~\ref{gradedmodule} that
$$v_p(\alpha') = \dim_{\mathbb{F}_p}\Big( \bigoplus_{d < a \leqslant e}
A(u)_{a,i} \, \Big).$$
The $\mathbb{F}_p$-vector space on the right-hand side has a basis
given by the elements $x_r^k \sigma_r^l$, where $1 \leqslant r
\leqslant u$, and where $k$ and $l$ are non-negative integers such
that $k + l = i$ and $d < p^{r-1}k \leqslant e$. Therefore, we have 
$$\begin{aligned}
{} & v_p(\alpha')
{} = \sum_{1 \leqslant r \leqslant u} \operatorname{card} \{ 0
\leqslant k \leqslant  i \mid d < p^{r-1}k \leqslant e \} 
 = \sum_{1 \leqslant r \leqslant u} \operatorname{card} \{ 1 \leqslant
 k \leqslant i \mid d < p^{r-1}k \leqslant e \} \cr
{} & = \sum_{1 \leqslant k \leqslant i} \operatorname{card}\{1
\leqslant r \leqslant u \mid d < p^{r-1}k \leqslant e \} 
 = \sum_{0 \leqslant h < i} \operatorname{card} \{1 \leqslant r
 \leqslant u \mid d < p^{r-1}(h+1) \leqslant e \} \cr
{} & = \sum_{0 \leqslant h < i} \big( \operatorname{card} \{1
\leqslant r \leqslant u \mid d < p^{r-1}(h+1) \} - \operatorname{card}
\{1 \leqslant r \leqslant u \mid e < p^{r-1}(h+1) \} \big), \cr
\end{aligned}$$
where, we recall, $d = d_p(m,u,j)$ and $e = d_p(n,u,j)$. The
inequality $d < p^{r-1}(h+1)$ is equivalent to the inequality
$(p^{u-1}j - 1)/m < p^{r-1}(h+1)$ which, in turn, is equivalent to the
inequality
$$p^{u-r}j < m(h+1).$$
Suppose that $s = s_p(m,h,j)$ satisfies $1 \leqslant s \leqslant
u$. Then $p^{s-1}j \leqslant m(h+1) < p^sj$. Therefore, the inequality 
$p^{u-1}j < m(h+1)$ is equivalent to the inequality $u-r \leqslant s -
1$. Hence,
$$\operatorname{card} \{ 1 \leqslant r \leqslant u \mid d <
p^{r-1}(h+1) \} = \operatorname{card} \{ u - (s - 1) \leqslant r
\leqslant u\} = s = s_p(m,h,j).$$ 
Finally, if $m(h+1) < j$, we also find that
$$\operatorname{card} \{ 1 \leqslant r \leqslant u \mid d <
p^{r-1}(h+1) \} = 0 = s_p(m,h,j).$$
The proof that $\operatorname{card} \{ 1 \leqslant r \leqslant u \mid 
e < p^{r-1}(h+1) \} = s_p(n,h,j)$ is similar.
\end{proof}

\begin{proof}[of Thm.~\ref{main}]We prove the equivalent
Thm.~\ref{ptypicalmain}. Suppose that $A = \mathbb{F}_p$. We first
identify the sequences~(\ref{m}) and~(\ref{mdrw}) and the
sequences~(\ref{n}) and~(\ref{ndrw}). Let $m \geqslant 1$ and
$i \geqslant 0$ be integers, and let $j \geqslant 1$ be an integer not
divisible by $p$. We know from Cor.~\ref{limitreachedcor} that, if $u'
\geqslant 1$ is an integer such that $p^{u'}j > m(i+1)$, then the
canonical projection 
$$\operatorname{pr}_{u'} \colon \lim_R
\operatorname{TR}_{2i-\lambda_d}^u(\mathbb{F}_p;p) \to 
\operatorname{TR}_{2i-\lambda_{d'}}^{u'}(\mathbb{F}_p;p),$$
where $d = d_p(m,u,j)$ and $d' = d_p(m,u',j)$, is an isomorphism.
Hence, the composition of the isomorphism $\sigma_p(d',u',i)$ of
Prop.~\ref{TRperfect} and the inverse of the isomorphism
$\operatorname{pr}_{u'}$ defines an isomorphism
$$\tau_p(m,i,j,u') \colon W_s(\mathbb{F}_p) \xrightarrow{\sim}
\lim_R \operatorname{TR}_{2i-\lambda_d}^u(\mathbb{F}_p;p),$$
where $s = s_p(m,i,j)$. Then, if $p^{u'} > m(i+1)$, we define
$$\tau_p(m,i,u') \colon \bigoplus_{j \in I_p} W_s(\mathbb{F}_p) \to
\prod_{j \in I_p} \lim_R
\operatorname{TR}_{2i-\lambda_d}^u(\mathbb{F}_p;p)$$ 
to be the isomorphism that takes the $j$th summand on the left-hand
side to the $j$th factor on the right-hand side by the map
$\tau_p(m,i,j,u')$. Similarly, we define
$$\tau_p'(m,i,u') \colon \bigoplus_{j \in m'I_p} W_{s-v}(\mathbb{F}_p)
\to \prod_{j \in m'I_p} \lim_R
\operatorname{TR}_{2i-\lambda_d}^{u-v}(\mathbb{F}_p;p),$$ 
where $m = p^vm'$ with $m'$ not divisible by $p$, where $s =
s_p(m,i,j)$, and where $d = d_p(m,u,j)$, to be the isomorphism that
takes the $j$th summand on the left-hand  side to the $j$th factor on
the right-hand side by the isomorphism given by the composition of
$\sigma_p(d',u'-v,i)$ and the inverse of the isomorphism
$$\operatorname{pr}_{u'} \colon \lim_R
\operatorname{TR}_{2i-\lambda_d}^{u-v}(\mathbb{F}_p;p) \to 
\operatorname{TR}_{2i-\lambda_{d'}}^{u'-v}(\mathbb{F}_p;p).$$
Then, for every $u'$ such that $p^{u'} > m(i+1)$, the diagram
$$\xymatrix{
{ \displaystyle{\bigoplus_{j \in m'I_p} W_{s-v}(\mathbb{F}_p) } }
\ar[r]<1ex>^(.52){m'V^v} \ar[d]^{\tau_p'(m,i,u')} &
{ \displaystyle{\bigoplus_{j \in I_p} W_s(\mathbb{F}_p) } }
\ar[d]^{\tau_p(m,i,u')} \cr
{ \displaystyle{ \prod_{j \in m'I_p} \lim_R
    \operatorname{TR}_{2i-\lambda_d}^{u-v}(\mathbb{F}_p;p) } } 
\ar[r]<1ex>^(.52){m'V^v} &
{ \displaystyle{ \prod_{j \in I_p} \lim_R
   \operatorname{TR}_{2i-\lambda_d}^u(\mathbb{F}_p;p) } } \cr 
}$$
commutes. This identifies the sequences~(\ref{m})
and~(\ref{mdrw}). However, we do not know that the family of
isomorphisms $\sigma_p(a,u,i)$ in Prop.~\ref{TRperfect} can be chosen
with the additional property that, if both $p^{u'} > m(i+1)$ and
$p^{u''} > m(i+1)$, the isomorphisms $\tau_p(m,i,u')$ and
$\tau_p(m,i,u'')$ are equal. Therefore, we choose $u' = u'(m,i)$ to be
the unique integer that satisfies
$$p^{u'-1} \leqslant m(i+1) < p^{u'}$$
and use the isomorphisms $\tau_p(m,i,u')$ and $\tau_p'(m,i,u')$ to
identify~(\ref{m}) and~(\ref{mdrw}). We use the same $u' = u'(m,i)$
and the isomorphisms $\tau_p(n,i,u')$ and $\tau_p'(n,i,u')$ to
identify the sequences~(\ref{n}) and~(\ref{ndrw}). In particular, the
sequences~(\ref{mdrw}) and~(\ref{ndrw}) are exact as proved
in~\cite[Thm.~4.2.10]{hm1}. Finally, Prop.~\ref{finalformula} shows
that there is a commutative diagram
$$\xymatrix{
{ W_s(\mathbb{F}_p) } \ar[r]^{R^{s-t}} \ar[d]^{\tau_p(m,i,j,u')} &
{ W_t(\mathbb{F}_p) } \ar[r]^{m_{\alpha}} &
{ W_t(\mathbb{F}_p) } \ar[d]^{\tau_p(n,i,j,u')} \cr
{ \lim_R \operatorname{TR}_{2i-\lambda_d}^u(\mathbb{F}_p;p) }
\ar[rr]^{\iota_p(m,n,i,j)} && 
{ \lim_R \operatorname{TR}_{2i-\lambda_e}^u(\mathbb{F}_p;p), } \cr
}$$
where the map $m_{\alpha}$ is given by multiplication by an element
$\alpha = \alpha_p(m,n,i,j) \in W(\mathbb{F}_p)$ with $p$-adic valution
$$v_p(\alpha) = \sum_{0 \leqslant h < i} \big(
s_p(m,h,j) - s_p(n,h,j) \big).$$
This completes the proof of Thm.~\ref{ptypicalmain} for
$A = \mathbb{F}_p$. 

Suppose next that $A$ is any regular noetherian ring that is also an 
$\mathbb{F}_p$-algebra. The groups
$\lim_R\operatorname{TR}_{q-\lambda_d}^u(A;p)$ were evaluated
in~\cite[Thm.~2.2.2]{hm2}, but see also~\cite[Thm.~16]{h3}. By the
universal property of the de~Rham-Witt complex, there is a canonical
map of graded rings
$$\xi_u \colon W_u\Omega_A^* \to \operatorname{TR}_*^u(A;p)$$
that commutes with $R$, $F$, $V$, and $d$. Let $m \geqslant 1$ and $i 
\geqslant 0$ be integers, and let $u' = u'(m,i)$. Let $j$ be an
integer not divisible by $p$, and let $s = s_p(m,i,j)$. We consider
the map
$$\tilde{\omega}_p(m,q,i,j,u') \colon
\lim_RW_u\Omega_A^{q-2i} \otimes_{W(\mathbb{F}_p)} W_s(\mathbb{F}_p) \to
\lim_R\operatorname{TR}_{q-\lambda_d}^u(A;p)$$
defined by the composition
$$\begin{aligned}
{} & \lim_RW_u\Omega_A^{q-2i} \otimes_{W(\mathbb{F}_p)}
W_s(\mathbb{F}_p) \to 
\lim_R\operatorname{TR}_{q-2i}^u(A;p) \otimes_{W(\mathbb{F}_p)}
\lim_R\operatorname{TR}_{2i-\lambda_d}^u(\mathbb{F}_p;p) \cr
{} & \to \lim_R\operatorname{TR}_{q-2i}^u(A;p) \otimes_{W(\mathbb{F}_p)}
\lim_R\operatorname{TR}_{2i-\lambda_d}^u(A;p) \to
\lim_R\operatorname{TR}_{q-\lambda_d}^u(A;p), \cr
\end{aligned}$$
where the first map is induced by the maps $\xi_u \otimes
\tau_p(m,i,j,u')$, where the second map is induced from the unit map 
$\eta \colon \mathbb{F}_p \to A$, and where the last map is induced by
the $T(A)$-module spectrum structure on $S^{\lambda_d} \wedge T(A)$. It
follows from~\cite[Thm.~16]{h3} that the maps
$\tilde{\omega}_p(m,q,i,j,u')$ factor through the canonical
projections
$$\operatorname{pr}_s \otimes \operatorname{id} \colon
\lim_RW_u\Omega_A^{q-2i} \otimes_{W(\mathbb{F}_p)} W_s(\mathbb{F}_p) \to
W_s\Omega_A^{q-2i}$$
and that the induced maps $\omega_p(m,q,i,j,u')$ define an isomorphism
$$\omega_p(m,q,j) = \bigoplus \omega_p(m,q,i,j,u') \colon
\bigoplus_{i \geqslant 0} W_s\Omega_A^{q-2i} \xrightarrow{\sim}
\lim_R\operatorname{TR}_{q-\lambda_d}^u(A;p),$$
where $u' = u'(m,i)$ and $s = s_p(m,i,j)$. We define
$$\omega_p(m,q) \colon \bigoplus_{i \geqslant 0}\bigoplus_{j \in I_p}
W_s\Omega_A^{q-2i} \to
\prod_{j \in I_p} \lim_R \operatorname{TR}_{q-\lambda_d}^u(A;p)$$ 
to be the isomorphism that takes the $j$ summand on the left-hand side
to the $j$th factor on the right-hand side by the isomorphism
$\omega_p(m,q,j)$. We define the isomorphism
$$\omega_p'(m,q) \colon \bigoplus_{i \geqslant 0}\bigoplus_{j \in m'I_p}
W_{s-v}\Omega_A^{q-2i} \to
\prod_{j \in m'I_p} \lim_R \operatorname{TR}_{q-\lambda_d}^{u-v}(A;p)$$
in a completely similar manner substituting the isomorphisms
$\tau_p'(m,i,j,u')$ for the isomorphisms $\tau_p(m,i,j,u')$ in the
definition of $\omega_p(m,q)$. Then the diagram
$$\xymatrix{
{ \displaystyle{\bigoplus_{i \geqslant 0}\bigoplus_{j \in m'I_p}
    W_{s-v}\Omega_A^{q-2i} } }
\ar[r]<1ex>^(.52){m'V^v} \ar[d]^{\omega_p'(m,q)} &
{ \displaystyle{\bigoplus_{i \geqslant 0}\bigoplus_{j \in I_p}
    W_s\Omega_A^{q-2i} } }
\ar[d]^{\omega_p(m,q)} \cr
{ \displaystyle{ \prod_{j \in m'I_p} \lim_R
    \operatorname{TR}_{q-\lambda_d}^{u-v}(A;p) } } 
\ar[r]<1ex>^(.52){m'V^v} &
{ \displaystyle{ \prod_{j \in I_p} \lim_R
    \operatorname{TR}_{q-\lambda_d}^u(A;p) } } \cr
}$$
commutes. Hence, the isomorphisms $\omega_p(m,q)$ and $\omega_p'(m,q)$ 
identify the sequences~(\ref{m}) and~(\ref{mdrw}). In particular, the
sequence~(\ref{mdrw}) is exact. Similarly, the isomorphisms
$\omega_p(n,q)$ and $\omega_p'(n,q)$ identify the sequences~(\ref{n})
and~(\ref{ndrw}). Finally, we have a commutative diagram
$$\xymatrix{
{ W_s\Omega_A^{q-2i} } \ar[r]^{R^{s-t}} \ar[d]^{\omega_p(m,q,i,j,u')} &
{ W_t\Omega_A^{q-2i} } \ar[r]^{m_{\alpha}} &
{ W_t\Omega_A^{q-2i} } \ar[d]^{\omega_p(n,q,i,j,u')} \cr
{ \lim_R \operatorname{TR}_{q-\lambda_d}^u(A;p) }
\ar[rr]^{\iota_p(m,n,i,j)} && 
{ \lim_R \operatorname{TR}_{q-\lambda_e}^u(A;p), } \cr
}$$
where the top horizontal map $m_{\alpha}$ is multiplication by the
\emph{same} element $\alpha = \alpha_p(m,n,i,j)$ of $W(\mathbb{F}_p)$
as in the basic case $A = \mathbb{F}_p$. This completes the proof of
Thm.~\ref{ptypicalmain}.
\end{proof}

\section{The divisor $\operatorname{div}(\alpha_p(m,n,i))$}
\label{divisorsection}

In this section we examine the divisor
$\operatorname{div}(\alpha_p(m,n,i))$ which appears in the 
statement of Thm.~\ref{main} and prove the precise Thm.~\ref{divisor}
below. We then use this result to derive Thm.~\ref{maink} of the
introduction and Thm.~\ref{milnorpart} below. We recall the integer
function $s_p(m,i,j)$ defined in~(\ref{s_p(m,i,j)}).

\begin{lemma}\label{exponential}Suppose that the integer $t =
s_p(n,i,j)$ satisfies 
$$\frac{m-n}{mn} \cdot \frac{p^t-1}{p-1} \cdot j
  \geqslant 2t.$$
Then $v_p(\alpha_p(m,n,i,j)) \geqslant s_p(n,i,j)$.
\end{lemma}

\begin{proof}We recall that
$$v_p(\alpha_p(m,n,i,j)) = \sum_{0 \leqslant h < i} \big (
s_p(m,h,j) - s_p(n,h,j) \big).$$
In this sum, every summand is non-negative, and the summand indexed by
$h$ is positive if and only if there exists an integer $r \geqslant 1$
with
$$n(h+1) < p^{r-1}j \leqslant m(h+1).$$
We estimate the number of indices $0 \leqslant h < i$ that satisfy
this inequality, for some $r \geqslant 1$. 
For a given $r \geqslant 1$, the number of integers $h$ that satisfy
the inequality
$$n(h+1) < p^{r-1}j \leqslant m(h+1)$$
is at least
$$\frac{(m-n)p^{r-1}j}{mn} - 1.$$
If $t = s_p(n,i,j) = 0$, then the statement of the lemma is trivially
true. So assume that $t > 0$. Then $t$ is the unique integer that
satisfies $p^{t-1}j \leqslant n(i+1) < p^tj$. Now, if $1 \leqslant r
\leqslant t$, then
$$n(h+1) < p^{r-1}j \leqslant p^{t-1}j \leqslant n(i+1),$$
implies that $h < i$. Therefore, we find
$$v_p(\alpha_p(m,n,i,j)) \geqslant \sum_{1 \leqslant r \leqslant t}
\big( \frac{(m-n)p^{r-1}j}{mn} - 1 \big) 
= \frac{m-n}{mn} \cdot \frac{p^t - 1}{p-1} \cdot j - t.$$
Hence, if the latter integer is greater than or equal to $t$,
or equivalently, if
$$\frac{m-n}{mn} \cdot \frac{p^t - 1}{p-1} \cdot j
\geqslant 2t,$$ 
then $v_p(\alpha_p(m,n,i,j)) \geqslant s_p(n,i,j)$.
\end{proof}

\begin{lemma}\label{stayszero}Suppose that $i \geqslant n/(m-n)$ and
$v_p(\alpha_p(m,n,i-1,j)) \geqslant s_p(n,i-1,j)$. Then also
$v_p(\alpha_p(m,n,i,j)) \geqslant s_p(n,i,j)$.
\end{lemma}

\begin{proof}The assumption that $i \geqslant n/(m-n)$ implies that
$mi \geqslant n(i+1)$ which, in turn, implies that $s_p(m,i-1,j)
\geqslant s_p(n,i,j)$. Suppose $v_p(\alpha_p(m,n,i-1,j)) \geqslant
s_p(n,i-1,j)$. Then
$$v_p(\alpha_p(m,n,i,j)) = v_p(\alpha_p(s,n,i-1,j)) + s_p(m,i-1,j) -
s_p(n,i-1,j) \geqslant s_p(m,i-1,j).$$
The statement follows.
\end{proof}

\begin{theorem}\label{divisor}The divisor
$\operatorname{div}(\alpha_p(m,n,i))$ satisfies the following:

(i) For every pair of integers $m > n > 1$ and for every prime
number $p$, there exists an integer $i_0 = i_0(m,n,p)$ such that, for 
all integers $i \geqslant i_0$,
$$\operatorname{div}(\alpha_p(m,n,i)) \geqslant
\operatorname{div}(\mathbb{W}_{n(i+1)}(\mathbb{F}_p)).$$

(ii) For every pair of integers $m > n > 1$, the integer $i_0(m,n,p)$
tends to infinity as the prime number $p$ tends to infinity.

(iii) For every integer $n > 1$ and every prime number $p$, there
exists an integer $m_0 = m_0(n,p)$ such that, for all $m \geqslant
m_0(n,p)$ and for all integers $i > 0$,
$$\operatorname{div}(\alpha_p(m,n,i)) \geqslant
\operatorname{div}(\mathbb{W}_{n(i+1)}(\mathbb{F}_p)).$$
\end{theorem}

\begin{proof}We first prove the statement~(i), which, by
Lemma~\ref{ptypicaldivisor} is equivalent to the statement that there 
exists an integer $i_0 = i_0(m,n,p)$ such that, for all positive
integers $j$ not divisible by $p$, $v_p(\alpha_p(m,n,i,j)) \geqslant
s_p(n,i,j)$. Suppose first that $j \geqslant 2mn/(m-n)$. Then
$$\frac{m-n}{mn} \cdot \frac{p^t - 1}{p-1} \cdot j \geqslant 2 \cdot 
\frac{p^t-1}{p-1}$$
and, for every prime number $p$ and every integer $t \geqslant 0$, the
right-hand side is greater than or equal to $2t$. Therefore,
Lemma~\ref{exponential} shows that $v_p(\alpha_p(m,n,i,j)) \geqslant
s_p(n,i,j)$, for all $i \geqslant 0$. Suppose next that $1 \leqslant j
< 2mn/(m-n)$. Then the inequality of the statement of
Lemma~\ref{exponential} is satisfied, if $t$ is large enough. Since $t
= s_p(n,i,j)$ tends to infinity as $i$ tends to infinity, it follows
that there exists an integer $i_0(m,n,p,j)$ such that
$v_p(\alpha_p(m,n,i,j)) \geqslant s_p(n,i,j)$, for all $i \geqslant
i_0(m,n,p,j)$. This proves that part~(i) of the statement holds with
$i_0(m,n,p)$ equal to the maximum of the integers $i_0(m,n,p,j)$,
where $1 \leqslant j < 2mn/(m-n)$. 

We next prove~(ii). Suppose that $p >
m(i+1)$. Then, for every integer $0 \leqslant h \leqslant i$, and for 
every integer $1 \leqslant j \leqslant n$ not divisible by $p$, 
$s_p(m,i,j) = s_p(n,i,j) = 1$, and hence, $v_p(\alpha_p(m,n,i,j)) <
s_p(n,i,j)$. Therefore, if $p > m(i+1)$, then $i_0(m,n,p) > i$. It
follows that $i_0(m,n,p) \geqslant [(p-1)/m]$ which tends to infinity
as $p$ tends to infinity.

Finally, we prove~(iii). For fixed $i$ and $j$, the integer
$s_p(m,i,j)$ tends to infinity as $m$ tends to infinity. Hence, there
exists $m_0 = m_0(n,p,i,j)$ such that $v_p(\alpha_p(m,n,i,j))
\geqslant s_p(n,i,j)$, for $m \geqslant m_0(n,p,i,j)$. Assume that
$m_0(n,p,i,j)$ is chosen minimal with this property. Then an induction
argument based on Lemma~\ref{stayszero} shows that $m_0(n,p,i,j)
\leqslant m_0(n,p,i-1,j)$, if $i \geqslant n/(m-n)$. Moreover, if $j >
n(i+1)$, then $s_p(n,i,j) = 0$, and therefore, $m_0(n,p,i,j) = n$.
This shows that~(iii) holds with $m_0(n,p)$ equal to the maximum of
the finitely many integers $m_0(n,p,i,j)$, where $0 \leqslant i
\leqslant n/(m-n)$ and $1 \leqslant j \leqslant n(i+1)$.
\end{proof}

\begin{lemma}\label{derhamwittzero}Let $A$ be an
$\mathbb{F}_p$-algebra and suppose that $A$ is generated by $N$
elements as an algebra over the subring $A^p$ of $p$th powers. Then,
for every subset $S \subset \mathbb{N}$ stable under division, the big
de~Rham-Witt group $\mathbb{W}_S\Omega_A^q$ is zero, if $q > N + 1$. 
\end{lemma}

\begin{proof}It suffices to show that, for all $u \geqslant 1$, the
$p$-typical de Rham-Witt group $W_u\Omega_A^q$ is zero, if
$q > N$ + 1; compare Sect.~\ref{ptypical} above. We first show that
$\Omega_A^q$ is zero, if $q > N$. By assumption, there exists a
surjective ring homomorphism $f \colon A^p[x_1,\dots,x_N] \to A$ from
a polynomial algebra in finitely many variables over $A^p$. The
induces map
$$f_* \colon \Omega_{A^p[x_1,\dots,x_N]}^q \to \Omega_A^q$$
is again surjective. Moreover, every element of of the domain is a sum 
of elements of the form
$$\omega = \eta dx_{i_1} \dots dx_{i_s}$$
where $0 \leqslant s \leqslant N$, where $1 \leqslant i_1 < \dots <
i_s \leqslant N$, and where $\eta \in \Omega_{A^p}^{q-s}$. Now, we
claim that $f_*(\eta) = 0$ unless $s = q$. Indeed, the element
$\eta \in \Omega_{A^p}^{q-s}$ can be written as a sum of elements of
the form $b_0 db_1 \dots db_{q-s}$, where $b_0, \dots, b_{q-s} \in
A^p$, and
$$f_*(b_0 db_1 \dots db_{q-s}) = f(b_0) d f(b_1) \dots d f(b_{q-s}).$$
But $f(b_i) = a_i^p$, for some $a_i \in A$, and hence, $d f(b_i) =
d(a_i^p) = pa_i^{p-1}da_i = 0$. It follows that the image of $f_*$ is
zero, if $q > N$. This shows that $\Omega_A^q$ is zero, if $q > N$, as
stated.

Finally, we show by induction on $u \geqslant 1$ that $W_u\Omega_A^q$
is zero, if $q > N+1$. The case $u = 1$ holds, since the canonical map 
$\Omega_A^q \to W_1\Omega_A^q$ is an isomorphism. Finally, the
induction step follows from the exact sequence
$$\xymatrix{
{ \Omega_A^q \oplus \Omega_A^{q-1} } \ar[rr]^(.55){V^{u-1}+dV^{u-1}} &&
{ W_u\Omega_A^q } \ar[r]^(.45){R} &
{ W_{u-1}\Omega_A^q } \ar[r] &
{ 0 } \cr
}$$
which is proved in~\cite[Prop.~3.2.6]{hm4}. 
\end{proof}

\begin{proof}[of Thm.~\ref{maink}]Since $m > n+1$, we can choose
$n < k < m$ and write $f$ as the composition of the canonical
projections $g \colon A[x]/(x^m) \to A[x]/(x^k)$ and $h \colon
A[x]/(x^k) \to A[x]/(x^n)$. We consider the following maps of
long-exact sequences from Thm.~\ref{main}.
$$\xymatrix{
{ \cdots } \ar[r] &
{ \bigoplus_{i \geqslant 0} \mathbb{W}_{m(i+1)}\Omega_A^{q-2i} }
\ar[r]^(.49){\varepsilon} \ar[d]^{g_*'} &
{ K_{q+1}(A[x]/(x^m),(x)) } \ar[r]^(.51){\partial} \ar[d]^{g_*} &
{ \bigoplus_{i \geqslant 0} \mathbb{W}_{i+1} \Omega_A^{q-1-2i} }
\ar[r] \ar[d]^{0} &
{ \cdots } \cr
{ \cdots } \ar[r] &
{ \bigoplus_{i \geqslant 0} \mathbb{W}_{k(i+1)}\Omega_A^{q-2i} }
\ar[r]^(.49){\varepsilon} \ar[d]^{h_*'} &
{ K_{q+1}(A[x]/(x^k),(x)) } \ar[r]^(.51){\partial} \ar[d]^{h_*} &
{ \bigoplus_{i \geqslant 0} \mathbb{W}_{i+1} \Omega_A^{q-1-2i} }
\ar[r] \ar[d]^{0} &
{ \cdots } \cr
{ \cdots } \ar[r] &
{ \bigoplus_{i \geqslant 0} \mathbb{W}_{n(i+1)}\Omega_A^{q-2i} }
\ar[r]^(.49){\varepsilon} &
{ K_{q+1}(A[x]/(x^n),(x)) } \ar[r]^(.51){\partial} &
{ \bigoplus_{i \geqslant 0} \mathbb{W}_{i+1} \Omega_A^{q-1-2i} }
\ar[r] &
{ \cdots } \cr
}$$
It follows from Lemma~\ref{derhamwittzero} and Thm.~\ref{divisor}(i)
that there exists an integer $q_0$ such that both the maps $g_*'$ and
$h_*'$ are zero, for $q \geqslant q_0$. Finally, a diagram chase
based on the diagram above shows that the composite map $f_* = h_*
\circ g_*$ is zero, for $q \geqslant q_0$. 
\end{proof}

Suppose that $A$ is a regular noetherian ring and an
$\mathbb{F}_p$-algebra. We let
$$\varepsilon_i \colon \mathbb{W}_{n(i+1)}\Omega_A^{q-2i} \to
K_{q+1}(A[x]/(x^n),(x))$$
be the $i$th summand of the map $\varepsilon$ in the long-exact
sequence
$$\xymatrix{
{ \cdots } \ar[r] &
{ \bigoplus_{i \geqslant 0} \mathbb{W}_{i+1}\Omega_A^{q-2i} }
\ar[r]^{V_n} &
{ \bigoplus_{i \geqslant 0} \mathbb{W}_{n(i+1)}\Omega_A^{q-2i} }
\ar[r]^{\varepsilon} &
{ K_{q+1}(A[x]/(x^n),(x)) } \ar[r] &
{ \cdots. }
}$$
We show that the image of the map $\varepsilon_0$ has the following
interpretation. 

\begin{theorem}\label{milnorpart}Let $A$ be a regular noetherian ring
that is also an $\mathbb{F}_p$-algebra, and let $n$ be a positive
integer. Then
$$\begin{aligned}
{} \bigcap_{m > n} & \operatorname{im} \big( f_* \colon
K_q(A[x]/(x^m),(x)) \to 
K_q(A[x]/(x^n),(x)) \big) \cr
{} = \hskip5pt & \operatorname{im} \big( \varepsilon_0 \colon
\mathbb{W}_n\Omega_A^{q-1} \to K_q(A[x]/(x^n),(x)) \big) \cr
\end{aligned}$$
\end{theorem}

\begin{proof}Let $g \colon A[x]/(x^m) \to A[x]/(x^k)$ and $h \colon
A[x]/(x^k) \to A[x]/(x^n)$ be the canonical projections, with $m > k > 
n$. A diagram based on the diagram from the proof of Thm.~\ref{maink}
above shows that, for every integer $q \geqslant 0$, we have
inclusions
$$\begin{aligned}
{} & \textstyle{
\operatorname{im}\big(
\bigoplus_{i \geqslant 0} \mathbb{W}_{m(i+1)}\Omega_A^{q-2i}
\xrightarrow{h_*'g_*'} 
\bigoplus_{i \geqslant 0} \mathbb{W}_{n(i+1)}\Omega_A^{q-2i} 
\xrightarrow{\varepsilon}
K_{q+1}(A[x]/(x^n),(x)) \big) 
} \cr
{} & \subset \operatorname{im}\big( K_{q+1}(A[x]/(x^m),(x))
\xrightarrow{f_*} K_{q+1}(A[x]/(x^n),(x)) \big) \cr
{} & \subset \textstyle{ 
\operatorname{im}\big(
\bigoplus_{i \geqslant 0} \mathbb{W}_{k(i+1)}\Omega_A^{q-2i}
\xrightarrow{h_*'}  
\bigoplus_{i \geqslant 0} \mathbb{W}_{n(i+1)}\Omega_A^{q-2i}
\xrightarrow{\varepsilon}
K_{q+1}(A[x]/(x^n),(x)) \big)
}
\end{aligned}$$
The maps $h_*'g_*'$ and $h_*'$ map the summand indexed by $i = 0$ to
the summand indexed by $i = 0$ by the appropriate restriction map,
which is surjective. Thm.~\ref{divisor} shows that, on summands
indexed by $i > 1$, the maps $h_*'g_*'$ and $h_*'$ are zero, if $m$ is
sufficiently large. The statement follows.
\end{proof}

\begin{remark}Let $k$ be a perfect field of positive characteristic
$p$, and let $n$ be a positive integer. It follows from
Thms.~\ref{main} and~\ref{divisor} that, for $m \gg n$, the composite
map
$$\mathbb{W}_{n(i+1)}(k) \xrightarrow{\varepsilon}
K_{2i+1}(k[x]/(x^n),(x)) \xrightarrow{\delta}
K_{2i}(k[x]/(x^m),(x^n))$$ 
induces an isomorphism
$$\Delta_n' \colon \mathbb{W}_{n(i+1)}(k)/V_n\mathbb{W}_{i+1}(k)
\xrightarrow{\sim} K_{2i}(k[x]/(x^m),(x^n)).$$
We expect that, for $i = 1$, the map $\Delta_n'$ is equal to the map
$\Delta_n$ of Stienstra~\cite[Cor.~3.6]{stienstra}.
\end{remark}

\begin{remark}One generally expects that, for every scheme $X$, there
exists an Atiyah-Hirzebruch type spectral sequence
$$E_{s,t}^2 = H^{t-s}(X;\mathbb{Z}(t)) \Rightarrow
K_{s+t}(X)$$
from the motivic cohomology of $X$ to the algebraic $K$-theory of
$X$. Such a sequence has been constructed in the case where $X$ is a
smooth scheme over a field~\cite{suslin5}. However, in general, there
presently does not exist a definition of motivic cohomology that could
serve as the $E^2$-term of a spectral sequence of this form. However,
see~\cite{blochesnault,rulling}. The long-exact sequence relating
$K(A[x]/(x^n),(x))$ and the big de~Rham-Witt groups and
Thm.~\ref{milnorpart} suggests that, for $X =
\operatorname{Spec}(A[x]/(x^n))$, where $A$ is a regular noetherian
$\mathbb{F}_p$-algebra, the spectral sequence takes the form
$$E_{s,t}^1 =  \mathbb{W}_s\Omega_A^{t-(s-1)}
\oplus \mathbb{W}_{(s+1)n}\Omega_A^{t-s} \Rightarrow
K_{s+t}(A[x]/(x^n),(x))$$
with the $d^1$-differential induced by the Verschiebung operator
$$V_n \colon \mathbb{W}_s\Omega_A^{t-(s-1)} \to
\mathbb{W}_{sn}\Omega_A^{t-(s-1)}$$
and with all higher differentials zero. The $E^2$-term is of this
hypothetical spectral sequence then is our candidate for the value of
the motivic cohomology groups of $A[x]/(x^n)$ relative  to the ideal
$(x)$. We note that this value of the motivic cohomology groups is in
agreement with the Beilinson-Soul\'{e} vanishing conjecture that
$H^q(A[x]/(x^n),(x);\mathbb{Z}(t))$ is zero, for $q \leqslant 0$.
\end{remark}

\begin{acknowledgements}This paper was written in part during a visit 
to the National University of Singapore. The author would like to
thank the university and Jon Berrick in particular for their
hospitality and support.
\end{acknowledgements}

\providecommand{\bysame}{\leavevmode\hbox to3em{\hrulefill}\thinspace}
\providecommand{\MR}{\relax\ifhmode\unskip\space\fi MR }
\providecommand{\MRhref}[2]{%
  \href{http://www.ams.org/mathscinet-getitem?mr=#1}{#2}
}
\providecommand{\href}[2]{#2}

\affiliationone{
Department of Mathematics\\
Massachusetts Institute of Technology\\
77 Massachusetts Avenue\\
Cambridge, Massachusetts\\
USA\\
\email{larsh@math.mit.edu}
}
\affiliationtwo{
Graduate School of Mathematics\\
Nagoya University\\
Chikusa-ku\\
Nagoya, Japan 464-8602\\
Japan\\
\email{larsh@math.nagoya-u.ac.jp}
}


\begin{thebibliography}{10}

\bibitem{bass}
H.~Bass, \emph{Algebraic {$K$}-theory}, W. A. Benjamin, Inc., New
  York-Amsterdam, 1968.

\bibitem{blochesnault}
S.~Bloch and H.~Esnault, \emph{The additive dilogarithm}, Documenta Math.,
  Extra Volume: Kazuya Kato's Fiftieth Birthday (2003), 131--155.

\bibitem{ega4}
A.~Grothendieck, \emph{{\'E}l{\'e}ments de g{\'e}om{\'e}trie alg{\'e}brique.
  {IV}. {{\'E}}tude locale des sch{\'e}mas et des morphismes de sch{\'e}mas
  ({Q}uatri{\`e}me {P}artie).}, Inst. Hautes {\'{E}}tudes Sci. Publ. Math.
  \textbf{32} (1967).

\bibitem{h}
L.~Hesselholt, \emph{On the $p$-typical curves in {Q}uillen's {$K$}-theory},
  Acta Math. \textbf{177} (1997), 1--53.

\bibitem{h3}
\bysame, \emph{{$K$}-theory of truncated polynomial algebras}, Handbook of
  {$K$}-theory, vol.~1, Springer-Verlag, New York, 2005, pp.~71--110.

\bibitem{hm1}
L.~Hesselholt and I.~Madsen, \emph{Cyclic polytopes and the {$K$}-theory of
  truncated polynomial algebras}, Invent. Math. \textbf{130} (1997), 73--97.

\bibitem{hm}
\bysame, \emph{On the {$K$}-theory of finite algebras over {W}itt vectors of
  perfect fields}, Topology \textbf{36} (1997), 29--102.

\bibitem{hm2}
\bysame, \emph{On the {$K$}-theory of nilpotent endomorphisms}, Homotopy
  methods in algebraic topology (Boulder, CO, 1999), Contemp. Math., vol. 271,
  Amer. Math. Soc., Providence, RI, 2001, pp.~127--140.

\bibitem{hm4}
\bysame, \emph{On the {$K$}-theory of local fields}, Ann. of Math. \textbf{158}
  (2003), 1--113.

\bibitem{hm3}
\bysame, \emph{On the de~{R}ham-{W}itt complex in mixed characteristic}, Ann.
  Sci. {\'E}cole Norm. Sup. \textbf{37} (2004), 1--43.

\bibitem{illusie}
L.~Illusie, \emph{Complexe de de {R}ham-{W}itt et cohomologie cristalline},
  Ann. Scient. \'Ec. Norm. Sup. (4) \textbf{12} (1979), 501--661.

\bibitem{illusieraynaud}
L.~Illusie and M.~Raynaud, \emph{Les suites spectrales associ{\'{e}}es au
  complexe de de~{R}ham-{W}itt}, Inst. Hautes {\'{E}}tudes Sci. Publ. Math.
  \textbf{57} (1983), 73--212.

\bibitem{mandellmay}
M.~A. Mandell and J.~P. May, \emph{Equivariant orthogonal spectra and
  {$S$}-modules}, Mem. Amer. Math. Soc., vol. 159, Amer. Math. Soc.,
  Providence, RI, 2002.

\bibitem{popescu}
D.~Popescu, \emph{General {N\'e}ron desingularization}, Nagoya Math. J.
  \textbf{100} (1985), 97--126.

\bibitem{rulling}
K.~R{\"{u}}lling, \emph{The generalized de~{R}ham-{W}itt complex over a field
  is a complex of zero-cycles}, J. Algebraic Geom. \textbf{16} (2007),
  109--169.

\bibitem{stienstra}
J.~Stienstra, \emph{On the {$K_2$} and {$K_3$} of truncated polynomial rings},
  Algebraic {$K$}-theory (Evanston, 1980), Lecture Notes in Math., vol. 854,
  Springer-Verlag, New York, 1981, pp.~409--455.

\bibitem{suslin5}
A.~A. Suslin, \emph{On the {G}rayson spectral sequence}, Number theory,
  algebra, and algebraic geometry, Proceedings of the Steklov Institute of
  Mathematics, vol. 241, 2003, pp.~202--237.

\end{thebibliography}
\end{document}